\documentclass{amsart}
\usepackage{epsf}
\usepackage[all]{xy}
\newcommand{\el}{\par \mbox{} \par \vspace{-0.5\baselineskip}}

\newtheorem*{theorem}{Theorem}
\newtheorem*{lemma}{Lemma}
\newtheorem*{proposition}{Proposition}

\theoremstyle{definition}

\theoremstyle{remark}

\numberwithin{equation}{section}

\newcommand{\internalcomment}[1]{}

\newcommand{\Z}{\mathbf{Z}}

\newcommand{\ko}{\: , \;}

\newcommand{\ol}{\overline}

\renewcommand{\tilde}[1]{\widetilde{#1}}

\newcommand{\ra}{\rightarrow}

\newcommand{\da}{\downarrow}
\newcommand{\ua}{\uparrow}

\newcommand{\arr}[1]{\stackrel{#1}{\rightarrow}}
\newcommand{\longarr}[1]{\xrightarrow{#1}}
\newcommand{\longlarr}[1]{\xleftarrow{#1}}
\newcommand{\iso}{\stackrel{_\sim}{\rightarrow}}

\newcommand{\larr}[1]{\stackrel{#1}{\leftarrow}}

\newcommand{\opname}[1]{\operatorname{\mathsf{#1}}}

\renewcommand{\mod}{\opname{mod}\nolimits}
\newcommand{\Mod}{\opname{Mod}\nolimits}

\newcommand{\obj}{\opname{obj}\nolimits}
\newcommand{\id}{\mathbf{1}}

\newcommand{\ten}{\otimes}
\newcommand{\tp}[1]{^{\ten #1}}

\newcommand{\im}{\opname{im}\nolimits}
\renewcommand{\ker}{\opname{ker}\nolimits}

\renewcommand{\H}[1]{{H}^{#1}}

\newcommand{\Hs}{{H}^*}

\newcommand{\ca}{{\mathcal A}}
\newcommand{\cb}{{\mathcal B}}
\newcommand{\cc}{{\mathcal C}}
\newcommand{\cd}{{\mathcal D}}

\newcommand{\ch}{{\mathcal H}}

\newcommand{\cs}{{\mathcal S}}
\newcommand{\ct}{{\mathcal T}}

\newcommand{\eps}{\varepsilon}

\renewcommand{\phi}{\varphi}

\newcommand{\Hom}{\opname{Hom}}
\newcommand{\HOM}{\opname{Hom^\bullet}}
\newcommand{\Ext}{\opname{Ext}}
\newcommand{\Hominf}{\raisebox{0ex}[2ex][0ex]{$\overset{\,\infty}{
                              \raisebox{0ex}[1ex][0ex]{$\mathsf{Hom}$}
                                                                 }$}}
\newcommand{\sHominf}{\raisebox{0ex}[2ex][0ex]{$\scriptsize \overset{\scriptsize \,\infty}{\raisebox{0ex}[1ex][0ex]{$\scriptsize \mathsf{Hom}$}}$}}

\newcommand{\tensinf}{\overset{\infty}{\ten}}
\newcommand{\Tor}{\opname{Tor}}

\newcommand{\filt}{\opname{filt}}
\newcommand{\tria}{\opname{tria}}
\newcommand{\tw}{\opname{tw}}

\newcommand{\centeps}[1]{\begin{array}{c} \epsfbox{#1} \end{array}}

\begin{document}
\title{Introduction to $A$-infinity algebras and modules}

\author{Bernhard Keller}
\address{UFR de Math\'ematiques\\
   UMR 7586 du CNRS \\
   Case 7012\\
   Universit\'e Paris 7\\
   2, place Jussieu\\
   75251 Paris Cedex 05\\
   France }
\thanks{Partially supported by the European Network `Algebraic
Lie Representations', Contract ERB-FMRX-CT97-0100}

\email{
\begin{minipage}[t]{5cm}
keller@math.jussieu.fr \\
www.math.jussieu.fr/ $\tilde{ }$ keller
\end{minipage}
}

\date{May 23, 1999, last modified on January 10, 2001}

\dedicatory{Dedicated to H.~Keller on the occasion of his seventy fifth 
birthday}

\begin{abstract}
These are expanded notes of four introductory
talks on $A_\infty$-algebras, their modules and their derived
categories. 
\end{abstract}

\maketitle

\tableofcontents

\section{Introduction}

\subsection{These notes}
These are expanded notes of a minicourse of three lectures 
given at the Euroconference 
`Homological Invariants in Representation
Theory' in Ioannina, Greece, March 16 to 21, 1999,
and of a talk at the Instituto de Matem\'aticas,
UNAM, M\'exico, on April 28, 1999. They present
basic results on $A_\infty$-algebras, their
modules and their derived categories.

\subsection{History}
$A_\infty$-spaces and  $A_\infty$-algebras (= sha algebras = 
strongly homotopy associative algebras) were invented at 
the beginning of the sixties by J.~Stasheff \cite{Stasheff63a}
\cite{Stasheff63b} as a tool in the study of `group-like'
topological spaces. In the subsequent two decades, 
$A_\infty$-structures found applications and
developments \cite{May72} \cite{BoardmanVogt73} 
\cite{Adams78} in homotopy theory;
their use remained essentially confined to this subject
(cf. however \cite{Ovsienko85}, \cite{Skoeldberg97}).
This changed at the beginning of the
nineties when the relevance of $A_\infty$-structures
in algebra, geometry and mathematical physics became
more and more apparent (cf. e.g. \cite{GetzlerJones90},
\cite{Stasheff90}, \cite{McCleary99}).
Of special influence was M.~Kontsevich's talk
\cite{Kontsevich94} at the International Congress
in 1994: Inspired by K.~Fukaya's preprint \cite{Fukaya93}
Kontsevich gave a conjectural interpretation of
mirror symmetry as the `shadow' of an equivalence
between two triangulated categories associated
with $A_\infty$-categories. His conjecture was
proved in the case of elliptic curves
by A.~Polishchuk and E.~Zaslow \cite{PolishchukZaslow98}.

\subsection{Scope and sources}
In these notes, we aim at presenting some
basic theorems on $A_\infty$-algebras, their modules
and their derived categories. The notes may serve as 
an introduction to a (very) small part of Kontsevich's 
lectures~\cite{Kontsevich98}.
Besides these, our main sources were 
\cite{Proute84},
\cite{Stasheff63b}, \cite{Smirnov80}, \cite{Kadeishvili82}, 
\cite{Kadeishvili85}, \cite{Kadeishvili87}. 
Some of the results we state are not proved
in this form in the literature. Complete proofs
will be given in \cite{Lefevre2000}.

\subsection{Contents}

In section 2, we motivate the introduction of
$A_\infty$-algebras and modules by two basic problems
from homological algebra: 
\begin{itemize}
\item[1.] the reconstruction of a complex from its homology, 
\item[2.] the reconstruction of the category of iterated 
          selfextensions of a module from its extension algebra.
\end{itemize}
Then we briefly present the
topological origin of $A_\infty$-structures.
Section 3 is devoted to $A_\infty$-algebras and their
morphisms. The central result is the
theorem on the existence of minimal models.
In sections 4 and 5, we introduce the derived category
of an $A_\infty$-algebra and we present the natural
framework for the solution of problem 1.
In section 6, we sketch the formalism of standard functors
and arrive at the solution of problem 2.
Section 7 presents the category of twisted objects,
which is of importance because it is `computable'.

\subsection{Omissions} 
The links with Morita theory for derived categories \cite{Rickard89}
or dg categories \cite{Keller94} have not been made
explicit, cf. however \cite{Keller99b}.

The relevance of $A_\infty$-algebras to boxes and matrix problems
was discovered by S.~Ovsienko \cite{Ovsienko85} (cf. also \cite{Ovsienko98}).
We only give a hint of this important development
in example \ref{quasihered}.

The notions of $A_\infty$-equivalence and of $A_\infty$-enhanced
triangulated categories, crucial for \cite{Kontsevich98}, will
be treated in \cite{Lefevre2000}.

\bigskip
\noindent
{\bf Acknowledgments.} The author is grateful to N.~Marmaridis
and C.~M.~Ringel
for the organization of the Euroconference in Ioannina, Greece,
where the material of sections 2 to 6 was first presented.
He expresses his sincere thanks to J.~A.~de~la~Pe\~{n}a
and C.~Geiss for the hospitality he enjoyed
at the Instituto de Matem\'aticas,
UNAM, M\'exico, where the material of section 7 was first
presented and a first draft of these notes completed.

\section{Motivation and topological origin}

\subsection{Motivation} \label{motivation}
Let $k$ be a field. The following
two problems will guide us in our exploration of
the world of $A_\infty$-algebras and modules.

\el
{\bf Problem 1.} Let $A$ be an associative
$k$-algebra with $1$, let 
\[
M=(\ldots \ra M^p \arr{d^p} M^{p+1} \ra \ldots)
\]
be a complex of (unital right) $A$-modules and let 
\[
\Hs M = \ker d^* / \im d^*\ko *\in\Z \ko
\]
be its homology. It is a naive but natural question
to ask whether, up to quasi-isomorphism, 
we can reconstruct the complex $M$ from
its homology $\Hs M$.
This is of course impossible, except under very restrictive
hypotheses, for example when $A$ is an hereditary
algebra (i.e. an algebra such that submodules of
projective $A$-modules are projective, see \cite{CartanEilenberg56}). 
So we ask {\em what additional structure is needed if
we want to reconstruct $M$ from its homology}.
The answer is that {\em $\Hs M$ carries a unique
(up to isomorphism) $A_\infty$-module
structure over $A$ (viewed as an $A_\infty$-algebra)
which encodes exactly the additional
information needed for this task}.

\el
{\bf Problem 2.} Let $B$ be an associative $k$-algebra
with $1$, let $M_1, \ldots, M_n$ be
$B$-modules and let $\filt(M_i)$ denote the full subcategory of
the category of right $B$-modules whose objects
admit finite filtrations with subquotients among
the $M_i$. In other words, $\filt(M_i)$ is the
closure under extensions of the $M_i$. A natural
question to ask is whether $\filt(M_i)$ is
determined by the extension algebra
\[
\Ext^*_{B}\,(M,M)\ko \mbox{ where } M=\bigoplus M_i \ko
\]
together with its idempotents corresponding to the $M_i$.
Again, the answer is no. This time, it is no even
if we assume that $B$ is hereditary. Again,
we ask what additional structure on the extension
algebra is needed to reconstruct the category
of iterated extensions.
The answer is that {\em the extension algebra admits
an $A_\infty$-structure which encodes the additional 
information needed for this task.}

\el
The aim of these lectures is to explain the answers
to the two problems and to show how they fit into
a general theory.

\subsection{Topological origin} At the beginning
of the sixties, J.~Stasheff invented $A_\infty$-spaces
and $A_\infty$-algebras \cite{Stasheff63a}, \cite{Stasheff63b}
as a tool in the study of `group-like' spaces. Let us
consider the basic example:

Let $(X,*)$ be a topological space with a base point $*$
and let $\Omega X$ denote the space of based loops in
$X$: a point of $\Omega X$ is thus a continuous map
$f : S^1 \ra X$ taking the base point of the circle
to the base point $*$. We have the composition
map
\[
m_2 : \Omega X \times \Omega X \ra \Omega X
\]
sending a pair of loops $(f_1,f_2)$ to the loop $f_1 * f_2 = m_2(f_1, f_2)$
obtained by running through $f_1$ on the first half of
the circle and through $f_2$ on the second half
\[
f_1 * f_2 \quad \quad \begin{array}{c}\epsfbox{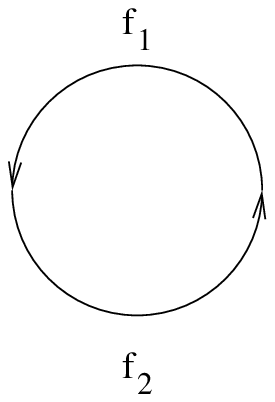}\end{array}
\]
This composition is non associative: for three loops
$f_1$, $f_2$, $f_3$, the composition $(f_1*f_2)*f_3$
runs through $f_1$ on the first {\em quarter} of the circle
whereas the composition $f_1*(f_2*f_3)$ runs through
$f_1$ on the first {\em half} of the circle. We symbolize the
two possibilities by the two binary trees with three
leaves
\begin{eqnarray*}
(f_1*f_2)*f_3 & \quad & \begin{array}{c} \epsfbox{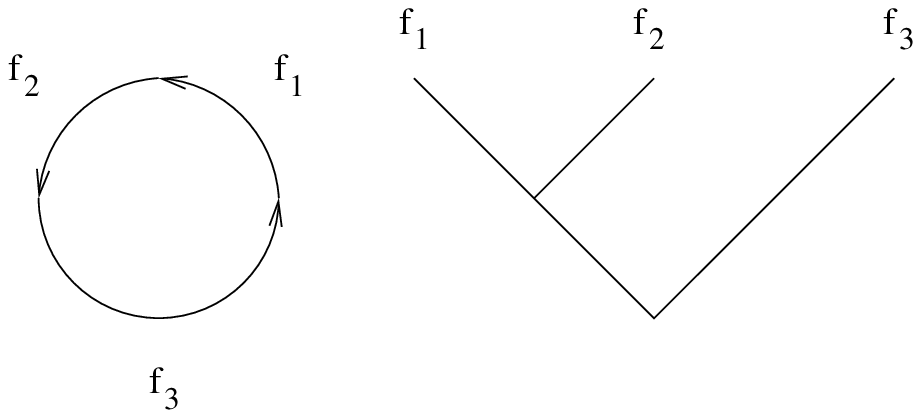} \end{array}\\
f_1*(f_2*f_3) & \quad & \begin{array}{c} \epsfbox{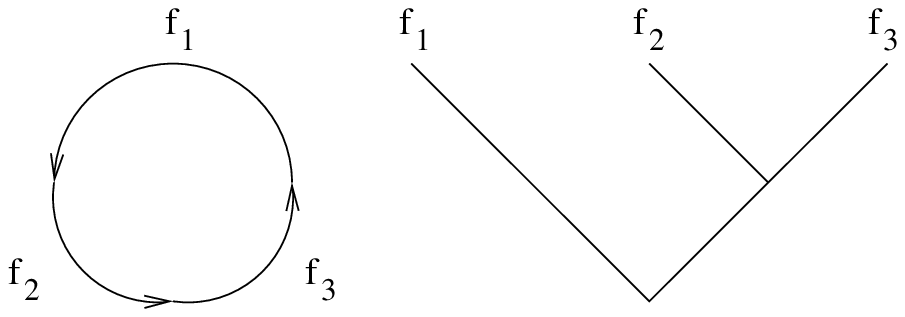} \end{array}
%% Insert ainf03.ps
\end{eqnarray*}
Clearly, there is an homotopy 
\[
m_3 : [0,1] \times \Omega X \times \Omega X \ra \Omega X.
\]
joining the two possibilities
of composing three loops. 
When we want to compose $4$ factors, there are $5$ 
possibilities corresponding to the $5$ binary trees with
$4$ leaves. Using $m_3$, we obtain two paths of homotopies
joining the composition $(f_1, f_2, f_3, f_4) \mapsto 
((f_1*f_2)*f_3)*f_4$ to the composition
$(f_1, f_2, f_3, f_4) \mapsto f_1*(f_2*(f_3*f_4))$.
\[
\epsfbox{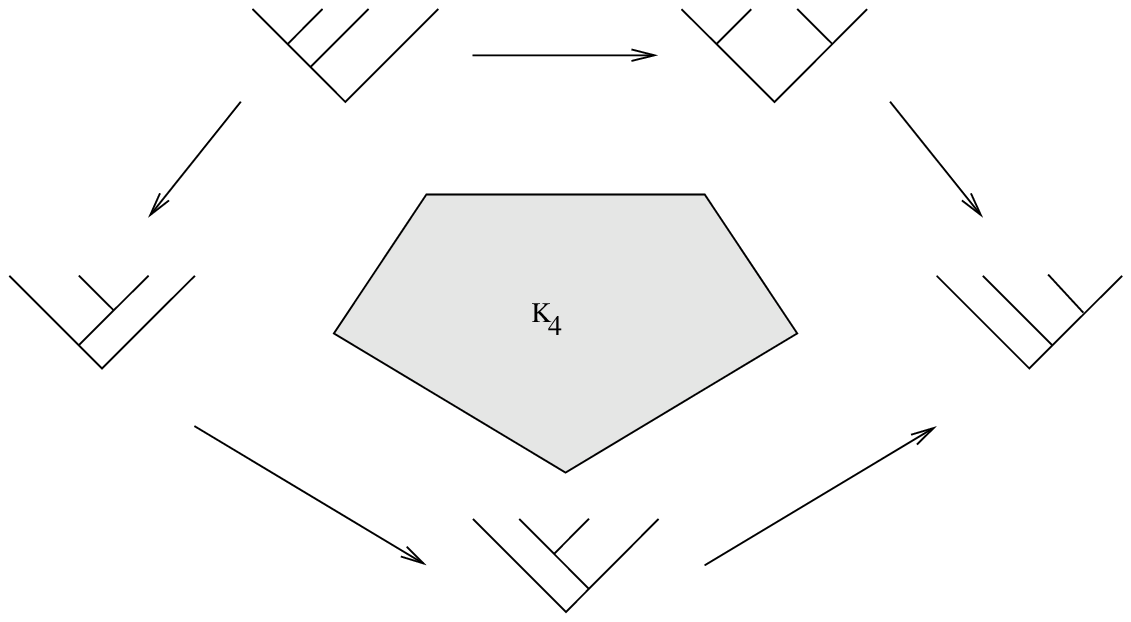}
\]
These two paths are clearly homotopic. Let us denote
an homotopy by 
\[
m_4 :K_4 \times (\Omega X)^4 \ra \Omega X
\]
where $K_4$ denotes the pentagon bounded by the
two paths. When we want to compose $5$ factors,
there are $14$ possibilities corresponding to
the $14$ binary trees with $5$ leaves. Using
$m_4$ and $m_3$, we obtain paths linking the compositions
and faces linking the paths. The resulting sphere is
the boundary of the polytope $K_5$ depicted 
\cite{Stasheff90} below. 
\[
\epsfbox{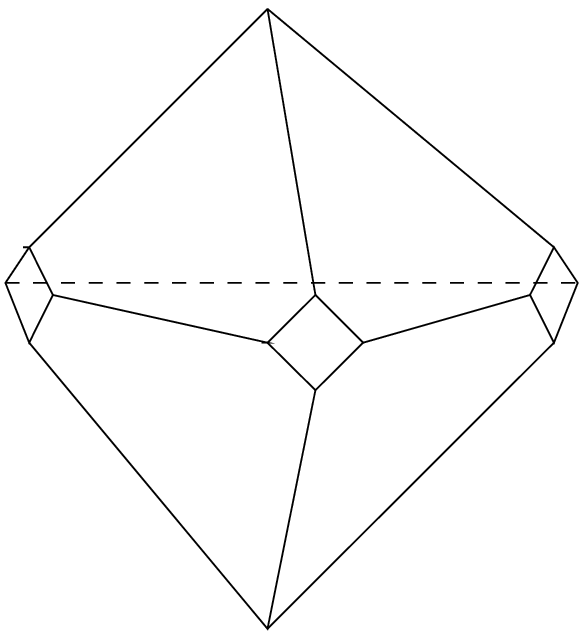}
\]
The pentagonal faces are copies of $K_4$. 
\internalcomment{
A rectangular face appears for example in the following
situation
\[
\epsfbox{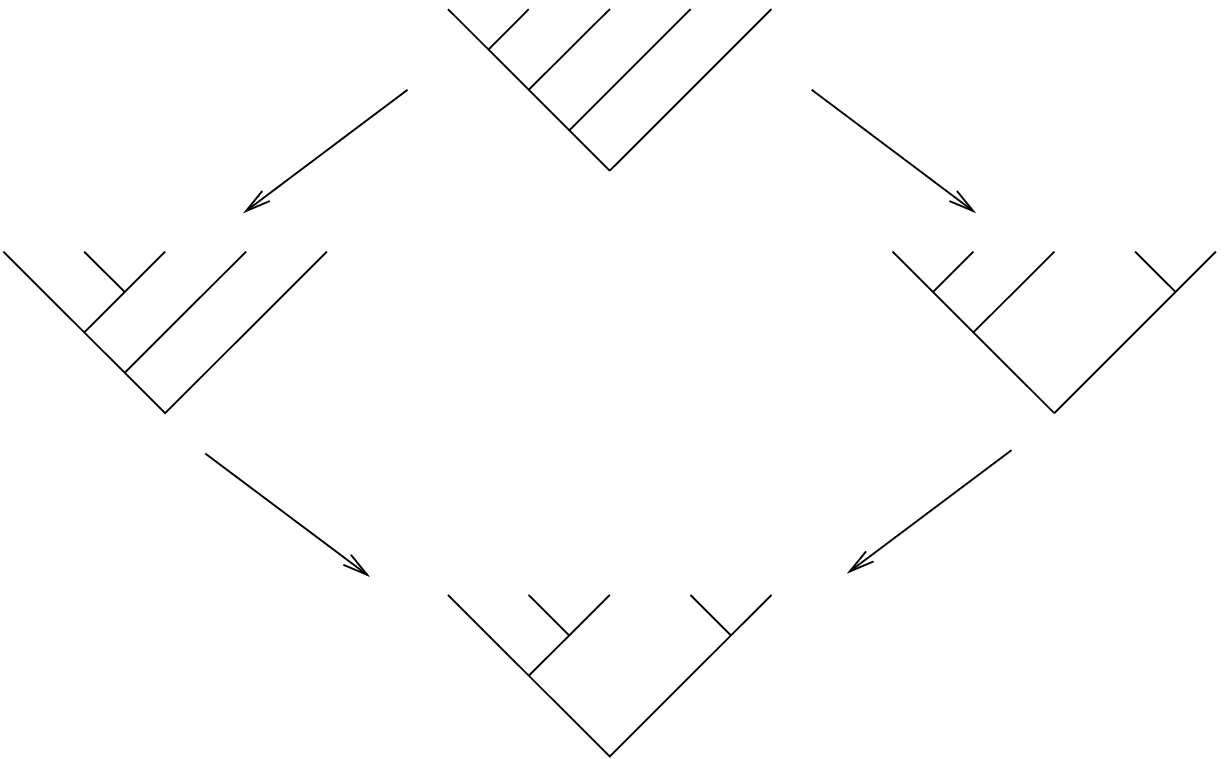}
\]
}

More generally, Stasheff \cite{Stasheff63a} 
defined polytopes $K_n$ of dimension $n-2$ for all $n\geq 2$ 
(we have $K_2=*$, $K_3=[0,1]$) 
and he defined an {\em $A_\infty$-space} to
be a topological space $Y$ endowed with maps
$m_n : K_n \times Y^n \ra Y$, $n\geq 2$,
satisfying suitable compatibility conditions
and admitting a `strict unit'.
The loop space $\Omega X$ is the prime example of 
such a space $Y$. Conversely \cite[2.2]{Adams78}, 
a topological space that admits the structure of
an $A_\infty$-space and whose connected components
form a group is homotopy equivalent
to a loop space.

If $Y$ is an $A_\infty$-space, the singular chain
complex of $Y$ is the paradigmatic example of
an $A_\infty$-algebra \cite{Stasheff63b}.

\internalcomment{
Let $C_*(Y)$ denote the singular chain complex on $Y$.
Then the map $m_n$ induces a map
\[
C_*(K_n \times Y^n) \ra C_*(Y).
\]
We compose this with the shuffle product map
(which is associative, cf. Loday, 1.6.11)
\[
C_*(K_n) \otimes C_*(Y)\tp{n} \ra C_*(K_n \times Y^n)
\]
to obtain the map
\[
C_*(K_n) \otimes C_*(Y)\tp{n} \ra C_*(Y).
\]
These maps satisfy certain compatibilities which
express that $C_*(Y)$ is an algebra over the 
asymmetric operad formed by the $C_*(K_n)$, $n\geq 1$.
The components $A_\infty(n)$ of the $A_\infty$-operad are not the $C_*(K_n)$. 
Rather, they are obtained as the cellular chain complexes
associated with the natural cell decomposition
of the $K_n$. By choosing a simplicial subdivision
of the cells, we obtain a system of
homotopy equivalences $A_\infty(n) \ra C_*(K_n)$
which is a quasi-isomorphism of operads.
This yields the $A_\infty$-algebra structure
on $C_*(Y)$.
}

\section{$A_\infty$-algebras and their morphisms}

\subsection{Definition of an $A_\infty$-algebra} 
\label{defainfalg} Let $k$ be a field.
An {\em $A_\infty$-algebra over $k$} (also called a 
`strongly homotopy associative algebra' or an `sha algebra')
is a $\Z$-graded vector space
\[
A=\bigoplus_{p\in\Z} A^p
\]
endowed with graded maps (=homogeneous $k$-linear maps)
\[
m_n : A\tp{n} \ra A \ko n\geq 1 \ko
\]
of degree $2-n$ satisfying the following
relations
\begin{itemize}
\item We have $m_1 \,m_1=0$, i.e. $(A,m_1)$ is a differential 
complex.
\item We have
\[
m_1 \, m_2 = m_2 \,(m_1 \ten \id + \id\ten m_1)
\]
as maps $A\tp{2} \ra A$. Here $\id$ denotes the identity
map of the space $A$. So $m_1$ is a (graded) derivation with 
respect to the multiplication $m_2$.
\item We have
\begin{multline*}
\quad \quad m_2 (\id\ten m_2 - m_2 \ten \id) \\
= m_1 m_3 + m_3\, (m_1 \ten \id \ten \id + 
                \id\ten m_1 \ten \id + \id \ten \id\ten m_1)
\end{multline*}
as maps $A\tp{3} \ra A$. Note that the left hand side is the 
associator for $m_2$ and that the right hand side
may be viewed as the boundary of $m_3$ in the morphism
complex $\HOM_k\,(A\tp{3},A)$ (the definition of the morphism
complex is recalled below). This implies that $m_2$ is 
associative up to homotopy. 
\item  More generally, for $n\geq 1$, we have
\[
\sum (-1)^{r+st} \,m_u \,(\id\tp{r} \ten m_s \ten \id\tp{t}) = 0
\]
where the sum runs over all decompositions $n=r+s+t$ and
we put $u=r+1+t$. 
\end{itemize}
We have adopted the sign conventions of Getzler-Jones
\cite{GetzlerJones90}.
Note that when these formulas are applied to elements,
additional signs appear because of the Koszul sign rule:
\[
(f\ten g)(x\ten y)= (-1)^{|g||x|} f(x)\ten g(y)
\]
where $f,g$ are graded maps, $x$,$y$ homogeneous elements
and the vertical bars denote the degree. For example, we have
\[
(m_1\ten \id+\id\ten m_1)(x\ten y)=
m_1(x)\ten y + (-1)^{|x|} x\ten m_1(y)
\]
so that $m_1\ten \id+\id\ten m_1$ is the usual differential
on the tensor product. 
The defining identities are pictorially represented 
\cite{Kontsevich98} by 
\medskip
\[
\sum \pm \centeps{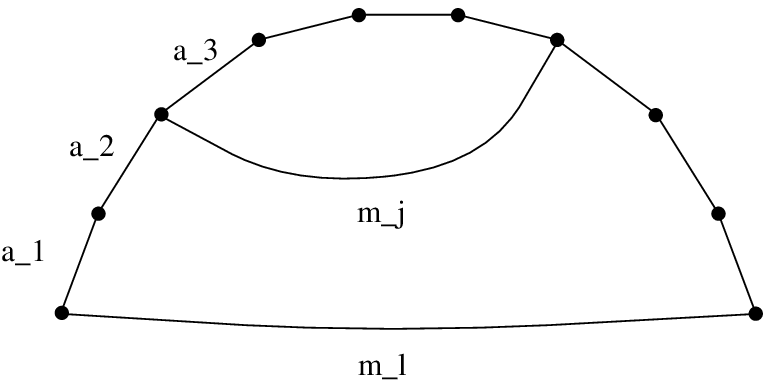}  = 0.
\]
We note three immediate consequences of the definition:
\begin{itemize}
\item[(1)] In general, an $A_\infty$-algebra $A$ is not associative but
its homology $\Hs A= \Hs(A,m_1)$ is an associative graded 
algebra for the multiplication induced by $m_2$. 

\item[(2)] If $A^p=0$ for all $p\neq 0$, then $A=A^0$ is an ordinary
associative algebra. Indeed, since $m_n$ is of degree $2-n$,
all $m_n$ other than $m_2$ have to vanish.

\item[(3)] If $m_n$ vanishes for all $n\geq 3$, then $A$ is an
associative dg (=differential $\Z$-graded) algebra and
conversely each dg algebra yields an $A_\infty$-algebra
with $m_n=0$ for all $n\geq 3$. 
\end{itemize}
Note that we have not defined unitality for $A_\infty$-algebras.
Cf. however \ref{AugmentedAlgebras} below.

\subsection{Link with deformation theory}
\label{deformationtheory}
The following class of examples, due to Penkava-Schwarz \cite{PenkavaSchwarz94},
shows the close link between deformation theory and $A_\infty$-algebras
(cf. also \cite{KontsevichSoibelman2000a}):
Let $B$ be an ordinary associative algebra, $n\geq 1$ an integer
and $\eps$ an indeterminate of degree $2-n$ and with $\eps^2=0$. 
An {\em infinitesimal $n$-deformation of $B$} is a 
$k[\eps]$-multilinear $A_\infty$-algebra structure on $B\ten k[\eps]$
which lifts the multiplication of $B$. Thus an infinitesimal $2$-deformation
is just an algebra deformation of $B$. It is well-known and easy to check
that infinitesimal $2$-deformations are in natural bijection with 
Hochschild $2$-cocycles $c: B\tp{2} \ra B$. Analogously, infinitesimal
$n$-deformations are in natural bijection with Hochschild $n$-cocycles.
Thus, by working in the class of $A_\infty$-algebras, we may interpret
Hochschild cocycles of {\em arbitrary degree} $\geq 1$ as deformations.

\subsection{Minimal models}
\label{minimalmodels}
Let us recall the dg algebra which is implicit in
problem 2: Let $B$ be
an ordinary associative $k$-algebra with $1$ and and $M$ a
$B$-module. Choose an injective resolution
\[
0 \ra M \ra I^0 \ra I^1 \ra \ldots
\]
of $M$ and denote by $I$ the complex
\[
0 \ra I^0 \ra I^1 \ra \ldots .
\]
Let $\HOM_B\,(I,I)$ be the morphism complex: recall that
its $n$th component is the group of graded $B$-linear maps $f: I \ra I$
of degree $n$ (no compatibility with $d_I$) and its 
differential is defined by 
$d(f) = d_I \circ f - (-1)^n f\circ d_I$.
Then $A=\HOM_B\,(I,I)$ endowed with its differential becomes
a dg algebra for the natural composition of graded maps and
its homology algebra is the extension algebra of $M$:
\[
\Hs A = \Ext^*_B\,(M,M).
\]
We claim that this has a canonical (up to equivalence) structure
of $A_\infty$-algebra with $m_1=0$ and $m_2$ the usual multiplication.
Indeed, more generally, we have the

\begin{theorem}[Kadeishvili \cite{Kadeishvili82}, see also
\cite{Kadeishvili80}, \cite{Smirnov80}, \cite{Proute84},
\cite{GugenheimLambeStasheff91}, \cite{JohanssonLambe96},
\cite{Merkulov98}]
\label{minimality}
Let $A$ be an $A_\infty$-algebra. Then the 
homology $\Hs A$ has an $A_\infty$-algebra
structure such that 
\begin{itemize}
\item[1)] we have $m_1=0$ and $m_2$ is induced by $m_2^A$ and
\item[2)] there is a quasi-isomorphism of $A_\infty$-algebras
$\Hs A \ra A$ lifting the identity of $\Hs A$.
\end{itemize}
Moreover, this structure is unique up to (non unique)
isomorphism of $A_\infty$-algebras.
\end{theorem}

The notion of morphism and quasi-isomorphism of $A_\infty$-algebras
will be defined below (\ref{morphisms}). 
The most explicit construction of the $A_\infty$-structure on $\Hs A$ 
can be found in Merkulov's \cite{Merkulov98}. His construction is explained
in terms of trees in \cite[6.4]{KontsevichSoibelman2000b}.

Let us introduce some terminology: An $A_\infty$-algebra with vanishing $m_1$
is called {\em minimal}.  In the situation of the theorem, the
$A_\infty$-algebra $\Hs A$ is called a {\em minimal model} for $A$.
The $A_\infty$-algebra $A$ is called {\em formal} if its minimal model
can be chosen to be an ordinary associative graded algebra, i.e. such
that $m_n=0$ for all $n\geq 3$.  A graded algebra $B$ is {\em
  intrinsically formal} if each $A_\infty$-algebra $A$ whose homology
is isomorphic to $B$ is formal. A sufficient condition for $B$ to be
intrinsically formal in terms of its Hochschild homology is given in
\cite{Kadeishvili88}, cf. also \cite[4.7]{SeidelThomas2000}.

It was shown in \cite{DeligneGriffithsMorganSullivan75}
that the dg algebra of smooth real forms on a compact K\"ahler
manifold is formal. In particular, this holds for
smooth, complex projective manifolds.

In contrast to the statement of the above theorem,
one can also show that each $A_\infty$-algebra $A$
admits an `anti-minimal' model, i.e. there is a dg
algebra $C$ and a quasi-isomorphism of $A_\infty$-algebras
$A \ra C$. Thus, passing from dg algebras to $A_\infty$-algebras
does not yield new quasi-isomorphism classes. What it does
yield is a {\em new description} of these classes 
by minimal models.

If, in the above example, we choose $B$ to be a 
finite-dimensional algebra and $M$ to be the sum of
the simple $B$-modules, one can show that
the dg algebra $A$ is formal iff $B$ is a Koszul algebra.
The smallest examples where $A$ is not formal is the 
algebra $B$ defined by the quiver with relations
(cf.~\cite{Ringel84}, \cite{GabrielRoiter92},
\cite{AuslanderReitenSmaloe95})
\[
\bullet \arr{\alpha} \bullet \arr{\beta} \bullet \arr{\gamma}\bullet
\ko \gamma\,\beta\,\alpha=0.
\]
Here the extension algebra $\Hs A$ is given by the quiver
with relations
\[
\centeps{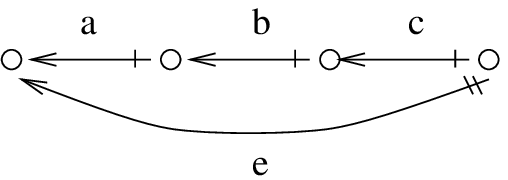} \ko ab=0\ko bc=0 \ko
\]
where the arrows $a,b,c$ are of degree $1$ and the arrow
$e$ is of degree $2$. 
We have $m_n=0$ for $n\neq 2,3$, $m_3(a,b,c)=e$ and
all other products vanish (the method for computing this
is given at the end of section \ref{morphisms}).
This example illustrates
the fact that the higher products of the $A_\infty$-algebra
$\Hs A$ lift the Massey products \cite{Stasheff63b} 
whenever these are defined.

\subsection{Morphisms of $A_\infty$-algebras}
\label{morphisms}
A {\em morphism of $A_\infty$-algebras} $f:A\ra B$ is a family
\[
f_n : A\tp{n} \ra B
\]
of graded maps of degree $1-n$ such that
\begin{itemize}
\item we have $f_1 m_1= m_1 f_1$, i.e. $f_1$ is a morphism
of complexes,
\item we have 
\[
f_1 m_2 = m_2 \, (f_1 \ten f_1) + m_1 f_2 + f_2 \, (m_1 \ten \id 
                               + \id\ten m_1) \ko
\]
which means that $f_1$ commutes with the multiplication
$m_2$ up to a homotopy given by $f_2$,
\item more generally, for $n\geq 1$, we have
\[
\sum (-1)^{r+st}\, f_u\,(\id\tp{r}\ten m_s \ten \id\tp{t}) 
=
\sum (-1)^s\, m_r\, (f_{i_1}\ten f_{i_2}\ten \ldots \ten f_{i_r})
\]
where the first sum runs over all decompositions
$n=r+s+t$, we put $u=r+1+t$, and the second sum runs over
all $1 \leq r\leq n$ and all decompositions $n=i_1+ \cdots +i_r$;
the sign on the right hand side is given by
\[
s=(r-1)(i_1 -1) + (r-2)(i_2-1)+ \cdots + 
2(i_{r-2}-1) + (i_{r-1}-1).
\]
\end{itemize}
The morphism $f$ is a {\em quasi-isomorphism} if $f_1$
is a quasi-isomorphism. It is {\em strict} if
$f_i=0$ for all $i\neq 1$. The {\em identity morphism}
is the strict morphism $f$ such that $f_1$ is the identity
of $A$.  The {\em composition} of two morphisms $f: B\ra C$ and 
$g:A \ra B$ is given by
\[
(f\circ g)_n = \sum (-1)^s f_r\, \circ \,(g_{i_1} \ten \cdots \ten
g_{i_r})
\]
where the sum and the sign are as in the defining identity.
Here is a (slightly stylized)
pictorial representation of the defining identity
for a morphism of $A_\infty$-algebras
\[
\sum \pm \centeps{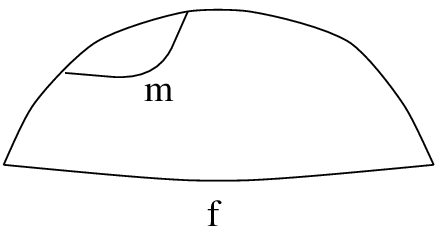} = \sum \pm \centeps{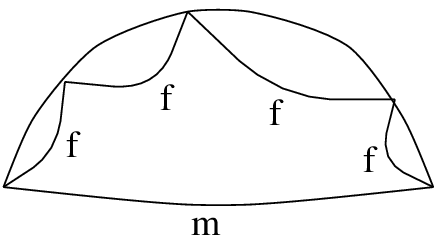}.
\]
$A_\infty$-algebras may be viewed as algebras over
an operad (cf. e.g. \cite{Markl99}). However, 
only strict morphisms of $A_\infty$-algebras may
then be viewed as morphisms of algebras.

As an illustration of the notion of morphism, suppose that $A$ is a dg
algebra and let us construct the multiplication $m_3$ on $\Hs A$ and,
at the same time, the first three terms of a quasi-isomorphism $f: \Hs
A \ra A$ lifting the identity of $\Hs A$.  For this, we view $\Hs A$
as a complex with zero differential and we choose a morphism of
complexes $f_1 :\Hs A \ra A$ inducing the identity in homology. Then
$f_1$ commutes with the multiplication up to boundaries, i.e. there
exists a graded map $f_2: \Hs A \ten \Hs A \ra A$ degree $-1$ such
that
\[
f_1 \,m_2 = m_2\,(f_1\ten f_1) + m_1\,f_2.
\]
Now we look for $f_3$ and $m_3$ such that
\begin{eqnarray*}
&   & f_1 \circ m_3 +f_2 \circ (m_2\ten\id - \id\ten m_2) \\
&   & \quad  +f_3 \circ (m_1\ten\id\tp{2}+\id\ten m_1 \ten\id 
                                       + \id\tp{2}\ten m_1)\\
& = & m_3 \circ (f_1\ten f_1 \ten f_1)
      +m_2\circ (f_1\ten f_2 - f_2\ten f_1)+m_1 f_3.
\end{eqnarray*}
In our situation, this simplifies into
\begin{eqnarray*}
f_1 \circ m_3 = m_2\circ (f_1\ten f_2 - f_2\ten f_1) + m_1 f_3.
\end{eqnarray*}
Now it is easy to check that the map 
$m_2\, (f_1 \ten f_2 - f_2 \ten f_1)$
has its image in the cycles of $A$. So there is indeed a map
$f_3 : (\Hs A)\tp{3} \ra A$ of degree $-2$ and a map
$m_3 : (\Hs A)\tp{3} \ra \Hs A$ of degree $-1$ such that
the identity holds.

\subsection{Augmented $A_\infty$-algebras}
\label{AugmentedAlgebras}
If an $A_\infty$-algebra $A$ has a unital homology algebra
$\Hs{A}$, one can often simplify the computation of the minimal
model for $A$ by passing to the category of augmented 
$A_\infty$-algebras. Here are the relevant definitions:

A {\em strictly unital $A_\infty$-algebra} is an $A_\infty$-algebra
$A$ endowed with an element $1_A$ of degree $0$ such that
we have $m_1(1_A)=0$, $m_2(1_A,a)=a=m_2(a,1_A)$ for all $a\in A$ and 
such that, for all $i>2$, and all $a_1, \ldots, a_i$ in $A$,
the product $m_i(a_1, \ldots, a_i)$ vanishes if one of the
$a_j$ equals $1_A$. If $A$ and $B$ are strictly unital
$A_\infty$-algebras, a {\em morphism} of $A_\infty$-algebras
$f: A \ra B$ is {\em strictly unital} if we have
$f_1(1_A)=1_B$ and if for all $i>1$ and all
$a_1, \ldots, a_i$, the element $f_i(a_1, \ldots, a_i)$ vanishes 
if one of the $a_j$ equals $1$. 
Each strictly unital $A_\infty$-algebra
is canonically endowed with a strict (hence strictly unital)
morphism $\eta : k \ra A$ mapping $1_k$ to $1_A$. 
It is {\em augmented} if it is moreover endowed with a strictly
unital morphism $\eps: A \ra k$ such that $\eps\circ \eta=\id_k$.
A {\em morphism of augmented $A_\infty$-algebras}
is a strictly unital morphism $f: A \ra B$ such that
$\eps_B \circ f = \eps_A$. It is then easy to check that the
functor $A \mapsto \ker \eps$ is an equivalence from
the category of augmented $A_\infty$-algebras
to the category of $A_\infty$-algebras. The quasi-inverse
functor sends an $A_\infty$-algebra $B$ to $k\oplus B$
endowed with the unique strictly unital structure
extending the structure of $B$.

Let $A$ be an augmented $A_\infty$-algebra and let
$\ol{A}=\ker\eps$. Clearly we obtain a minimal
model for $A$ by taking $k\oplus \ol{A}'$, where $\ol{A}'$
is a minimal model for $\ol{A}$ in the category
of $A_\infty$-algebras. The computation
of the minimal model for $\ol{A}$ is simpler than
that of $A$, because we do not have to deal with
the `loop' $1$. 

The construction of minimal models and the 
above discussion extend to the case where 
the ground ring $k$ is no longer a field but a semi-simple
ring $R$, for example a product of fields. We
then consider $A_\infty$-algebras in the monoidal
category of $R$-$R$-bimodules.

For example, suppose that $A$ is an ordinary algebra
augmented over a semi-simple ring $R$, say $A=R\oplus \ol{A}$.
Then the standard complex computing $\Ext_A(R,R)$ is of the
form
\[
0 \ra R \arr{0} \Hom_R(\ol{A}, \ol{A}) \ra \Hom_R(\ol{A}\ten_R \ol{A},
\ol{A}) \ra \ldots
\]
It has a canonical structure of dg algebra augmented over $R$.
Thus the $A_\infty$-algebra $\Ext_A(R,R)$ may be chosen to be
augmented over $R$ and, in particular, strictly unital.

Let us consider another example: Suppose that $k$ is a field and 
that $A$ is a unital dg $k$-algebra endowed with a decomposition
\[
1_A=e_1 + \ldots + e_N
\]
of $1_A$ into idempotents $e_i=e_i^2$ such that the
complexes $e_j A e_i$ are acyclic for $i<j$ and that
$e_i A e_i$ is quasi-isomorphic to $k\,e_i$. Clearly
$A$ is isomorphic to the matrix algebra
\[
\bigoplus_{i,j=1}^N e_j A e_i
\]
and this algebra is quasi-isomorphic to its subalgebra
of upper triangular matrices
\[
B=k e_1 \oplus \ldots \oplus k e_N \oplus \bigoplus_{i<j} e_j A e_i.
\]
Clearly $B$ is an augmented algebra over $R=k\times\ldots \times k$
($N$ factors). So we obtain a minimal model for $A$ as the
sum of $R$ with a minimal model for 
\[
\ol{B} = \bigoplus_{i<j} e_j A e_i.
\]
viewed as an $A_\infty$-algebra over $R$.
Let $\ol{B}'$ be this minimal model.
Then we have a decomposition 
\[
\ol{B}' = \bigoplus_{i<j} \Hs(e_j A e_i) = \bigoplus_{i<j} e_j \Hs(A) e_i
\]
of $R$-$R$-bimodules and the tensor power
\[
\ol{B}' \ten_R \ol{B}' \ten_R \ldots \ten_R \ol{B}'
\quad\mbox{($n$ factors)}
\]
is isomorphic to the sum of the
\[
e_{i_{n+1}} \Hs(A) e_{i_n} \ten_k e_{i_{n}} \Hs(A) e_{i_{n-1}} \ten
\ldots \ten_k e_{i_2} \Hs(A) e_{i_1}
\]
taken over all sequences $i_1 < i_2 < \ldots < i_{n+1}$. These spaces
vanish for $n\geq N$ so that in the minimal model we automatically
have $m_n=0$ for $n\geq N$. 
For example,
if $N=3$, i.e. if the unit of $A$ decomposes into $1=e_1 + e_2 + e_3$,
then $\ol{B}'$ is automatically formal.
Hence $\ol{B}$ and $A$ are also formal in this case.

The following lemma presents another case where the passage
to the category of augmented $A_\infty$-algebras is easy.
A {\em subalgebra} of an $A_\infty$-algebra $A$ is
a graded subspace $B$ such that $m_i$ maps $B\tp{i}$
to $B$ for all $i\geq 1$.

\begin{lemma} Let $A$ be a strictly unital $A_\infty$-algebra
(for example a unital dg algebra).
Suppose that $\H{i} A=0$ for $i<0$ and that $\H{0} A = k$.
Then there is a subalgebra $B\subset A$ whose inclusion
is a quasi-isomorphism and which is augmented.
\end{lemma}

%\begin{proof} We put $B^i=0$ for $i<0$, $B^0=k$, and 
%$B^i=A^i$ for $i\geq 2$. To define $B^1$, we choose
%decompositions of vector spaces
%\[
%A^1 = \ker m_1 \oplus V \quad\mbox{and}\quad
%\ker m_1 = W \oplus \im m_1.
%\]
%Note that $W$ is isomorphic to $\H{1}A$ as a $k$-vector
%space and that $m_1$ maps $V$ isomorphically onto its image.
%We put $B^1 = W \oplus V$. It is then clear that $(B, m_1)$ is
%a subcomplex of $(A,m_1)$ and that its inclusion is a
%quasi-isomorphism. It remains to be checked that
%$m_i$ maps $B\tp{i}$ to $B$ for all $i\geq 2$. 
%Suppose that $i=2$. Then $B^0\ten B^1$ and $B^1\ten B^0$ 
%are mapped to $A^1$ by $m_2$ and their images lie
%in $B^1$ since $1$ is a unit for $m_2$. If
%$j+k\geq 2$, then $B^j\ten B^k$ is mapped to
%$A^{j+k}=B^{j+k}$.
%Now suppose that $i>2$.
%Since $1$ is a strict unit for $A$, it is enough
%to show that $m_i$ maps $\ol{B}\tp{i}$ to $B$,
%where $\ol{B}$ is the sum of the $B^j$, $j>0$.
%Now $\ol{B}\tp{i}$ is concentrated in degrees
%$\geq i$. So its image under $m_i$ is concentrated
%in degrees $\geq 2$. But for $j\geq 2$, we have
%$B^j=A^j$.
%\end{proof}

\subsection{The bar construction} 
\label{barconstruction} Following \cite{Stasheff63b} and
\cite{Kadeishvili85} we will restate the definitions of
$A_\infty$-algebras and their morphisms in a more
efficient way.
Let $V$ be a $\Z$-graded vector space and let
\[
\ol{T}V= V \oplus V\tp{2} \oplus \cdots
\]
be the reduced tensor algebra (in the category of
graded spaces). We make $\ol{T}V$ into a graded
coalgebra via the comultiplication
\[
\Delta : \ol{T}V \ra \ol{T}V\ten\ol{T}V
\]
defined by
\[
\Delta(v_1, \ldots, v_n) = \sum_{1\leq i\leq n} 
(v_1, \ldots, v_i)\ten (v_{i+1}, \ldots, v_n)\; .
\]
Note that $\Delta(v)=0$ and $\Delta (v_1,v_2)=v_1\ten v_2$.
It is then easy to see that each graded map
$b: \ol{T} V \ra V$ of degree $1$ uniquely lifts
to a graded coderivation $b:\ol{T}V\ra\ol{T}V$
of degree 1, where {\em coderivation} means that
we have
\[
\Delta b = (b\ten \id + \id\ten b) \,\Delta.
\]
To see that this is plausible, one should remember that
algebra derivations $\ol{T}V \ra \ol{T}V$ are in
bijection with maps $V\ra \ol{T}V$. The component of
$b$ that maps $V\tp{n}$ to $V\tp{u}$ is
\[
\id\tp{r}\ten b_s \ten \id\tp{t}
\]
where $r+1+t=u$. Thus we have
a bijection between coderivations
$b:\ol{T}V \ra \ol{T}V$ of degree $1$ and families
of maps $b_n : V\tp{n}\ra V$, $n\geq 1$, of degree $1$.

Now let $A$ be another $\Z$-graded space and let
$SA$ be its suspension: $(SA)^p=A^{p+1}$. Define
a bijection between the families of maps
$b_n : (SA)\tp{n}\ra SA$ of degree $1$ and maps $m_n : A\tp{n}\ra A$,
$n\geq 1$, by the commutative square
\[
\begin{array}{rcl}
(SA)\tp{n} & \arr{b_n} & SA \\
s\tp{n}\ua &           & \ua s \\ 
A\tp{n}    & \arr{m_n} & A
\end{array}
\ko
\]
where $s: A \ra SA$ is the canonical map $a \mapsto a$
of degree $-1$. Note that $m_n$ is of degree $-n+1-(-1)=2-n$.
Then we have the following easy

\begin{lemma} The following are equivalent
\begin{itemize}
\item[(i)] The maps $m_n : A\tp{n} \ra A$ yield an
$A_\infty$-algebra structure on $A$.
\item[(ii)] The coderivation $b: \ol{T}SA \ra \ol{T}SA$ satisfies $b^2=0$.
\item[(iii)] For each $n\geq 1$, we have
\[
\sum b_u\,(\id\tp{r} \ten b_s \ten \id\tp{t}) =0 \ko
\]
where the sum runs over all decompositions $n=r+s+t$ and we
put $u=r+1+t$.
\end{itemize}
\end{lemma}

Note that no signs appear in the last formula. 
The lemma yields a bijection between $A_\infty$-algebra structures
on $A$ and coalgebra differentials on $\ol{T}SA$. 
It is also easy to check that if $A$ and $B$ are $A_\infty$-algebras,
then morphisms of $A_\infty$-algebras $f:A \ra B$ are in
bijection with morphisms of dg coalgebras $\ol{T}SA \ra \ol{T}SB$.

\subsection{Homotopy equivalences and quasi-isomorphisms}
\label{homotopyequivalences}
Let $A$ and $B$ be $A_\infty$-algebras and $f,g$ morphisms
from $A$ to $B$. Let $F,G$ denote the corresponding
morphisms of dg coalgebras $\ol{T}SA \ra \ol{T}SB$. 
One defines $f$ and $g$ to be {\em homotopic} if there
is a homogeneous map $H: \ol{T}SA \ra \ol{T}SB$ of
degree $-1$ such that 
\[
\Delta H = F\ten H + H \ten G \quad\mbox{and}\quad
F-G= b\circ H + H \circ b.
\]
It is easy to translate this into the existence
of a family of maps $h_n: A\tp{n} \ra B$ satisfying
suitable conditions. 

\begin{theorem}[Prout\'e \cite{Proute84}, cf. also \cite{Kadeishvili87}]
\begin{itemize}
\item[(a)] Homotopy is an equivalence relation on the
set of morphisms of $A_\infty$-algebras $A \ra B$.
\item[(b)] A morphism of $A_\infty$-algebras $s:A\ra B$
is a quasi-isomorphism iff it is an homotopy equivalence
(i.e. it becomes an isomorphism in the category obtained
by dividing the category of $A_\infty$-algebras by the
homotopy relation).
\end{itemize}
\end{theorem}

\section{$A_\infty$-modules and the derived category}
\label{ainfmodules}

\subsection{Reminder on the derived category} \label{reminder}
We refer to \cite[Ch. I]{KashiwaraSchapira90} for a 
concise presentation of derived categories and to
\cite{Weibel94} and \cite{GelfandManin96} for more leisurely
introductions. We refer to \cite{Spaltenstein88} for 
foundational material on resolutions of unbounded complexes.

Let $A$ be an associative $k$-algebra with $1$ and denote
by $\Mod A$ the category of unital right $A$-modules. The derived
category  $\cd A$ is classically 
constructed \cite{Verdier77} \cite{Verdier97} in three steps:

\bigskip
{\em 1st step:} We define the {\em category  $\cc A$} to be
the category of complexes
\[
\ldots \ra M^p \ra M^{p+1} \ra \ldots
\]
with morphisms of complexes as morphisms.

\bigskip
{\em 2nd step:} We define the {\em homotopy category $\ch A$}.
Its objects are the same as those of $\cc A$. The group
of morphisms $M \ra N$ in $\ch A$ is the group of morphisms
of complexes divided by the subgroup of null homotopic 
morphisms.

\bigskip
{\em 3rd step:} We define the {\em derived category $\cd A$}
to be the localization of the homotopy category $\ch A$ with
respect to all quasi-isomorphisms. Thus the objects of $\cd A$
are the same as those of $\ch A$ and its morphisms are
obtained from morphisms in the homotopy category by formally
inverting all quasi-isomorphisms.

To get an idea of what the morphisms in the derived category
are it is helpful to consider the following special case:
Let $M$ and $N$ be two $A$-modules. Identify $M$ with the
complex
\[
\ldots \ra 0 \ra M \ra 0 \ra \ldots
\]
concentrated in degree $0$. For any complex $K$ and any
$n\in\Z$, denote by $K[n]$ the complex $K$ shifted by
$n$ degrees to the left: $K[n]^p=K^{n+p}$ and 
$d_{K[n]}=(-1)^n d_K$.
Then there is a canonical quasi-isomorphism
\[
\Hom_{\cd A}\,(M,N[n]) = \Ext^n_A\,(M,N)
\]
valid for any $n\in \Z$ if we take the right hand side to
vanish for negative $n$.

\subsection{Generalization to $A_\infty$-algebras} 
\label{generalization}
Let $A$ be an $A_\infty$-algebra. 
We generalize the above step by step.

\bigskip
{\em 1st step: } An {\em $A_\infty$-module over $A$} (cf.
\cite{Kadeishvili86}) is a $\Z$-graded space $M$ endowed with maps
\[
m^M_n : M\ten A\tp{n-1} \ra M\ko n\geq 1\ko
\]
of degree $2-n$ such that  an identity of the same form
as the one in the definition
(\ref{defainfalg}) of an $A_\infty$-algebra holds. 
However, the term
\[
m_u\, (\id\tp{r} \ten m_s\ten \id\tp{t})
\]
has to be interpreted as
\[
m^M_u\, (\id\tp{r} \ten m_s\ten \id\tp{t})
\]
if $r>0$ and as
\[
m^M_u\, (m^M_s\ten \id\tp{t})
\]
if $r=0$.
Equivalently, the datum of an $A_\infty$-structure on a graded
space $M$ is the datum of a differential $m_1^M$ on $M$ and of
a morphism of $A_\infty$-algebras from $A$ to the opposite
of the dg algebra $\HOM_k(M,M)$ (cf. \ref{minimalmodels}).
\internalcomment{
A structure of $A_\infty$-module $b: TSA \ten SM \ra SM$ yields
a differential $b_1: SM \ra SM$, and hence $d_M = -b_1 : M \ra M$,
and a graded morphism
\[
\tau: \ol{T}SA \ra \HOM_k(M,M)
\]
of degree $+1$ such that 
\[
\tau\circ b = m_1 \circ \tau + m_2 * \tau \ko
\]
where $m_2 * \tau = \mu \circ (\tau\ten\tau)\circ \Delta$. 
The datum of $\tau$ is equivalent to that of a morphism
of $A_\infty$-algebras.
}
A {\em morphism of $A_\infty$-modules $f:L\ra M$} 
is a sequence of graded morphisms
\begin{equation} \label{MorphismAinfModule}
f_n : L \ten A\tp{n-1} \ra M 
\end{equation}
of degree $1-n$ such that for each $n\geq 1$, we have
\begin{equation} \label{CompatMorphismAinfModule}
\sum (-1)^{r+st}\,f_u \circ (\id\tp{r} \ten m_s \ten \id\tp{t})
=
\sum (-1)^{(r+1)s}\,m_u \circ (f_r\ten \id\tp{s}) \ko
\end{equation}
where the left hand sum is taken over all decompositions
$n=r+s+t$, $r,t\geq 0$, $s\geq 1$ and we put $u=r+1+t$;
and the right hand sum is taken over all decompositions
$n=r+s$, $r\geq 1$, $s\geq 0$ and we put $u=1+s$.
The morphism $f$ is a {\em quasi-isomorphism} if
$f_1$ is a quasi-isomorphism. The {\em identical morphism}
$f:M \ra M$ is given by $f_1=\id_M$, $f_i=0$, for all $i>0$.
The {\em composition $fg$} is given by
\[
(fg)_n=\sum (-1)^{(r-1)\,s} f_u (g_r \ten \id\tp{s})
\]
where the sum runs over all decompositions $n=r+s$ and
we put $u=1+s$.

In analogy with (\ref{barconstruction}), one can easily check
that $A_\infty$-structures on a graded space $M$ are in bijection
with comodule differentials on $SM\ten TSA$, where $TSA$ is
the coaugmented tensor coalgebra on $SA$ :
\[
TSA = k \oplus \ol{T}SA.
\]
The comultiplication 
on $TSA$ is given by
\[
\Delta(x) = 1\ten x + \Delta_{\ol{T}SA}(x) + x\ten 1.
\]
The differential of $TSA$ is induced by $b$ on
$\ol{T}SA$ and maps $k$ to zero. 
The comodule structure
on $SM\ten TSA$ is induced from the coalgebra structure
on $TSA$. Morphisms of $A_\infty$-modules then correspond
bijectively to morphisms of dg comodules.

%An $A_\infty$-module
%$M$ is \emph{$H$-unital} if $SM\ten TSA$ is acyclic.
%This holds for example if $\Hs M$ is $H$-unital
%in the sense of Wodzicki \cite{Wodzicki89} as
%a module over $\Hs A$. In particular, each complex of 
%unital modules over an ordinary associative algebra $A$ with 
%$1$ is an $H$-unital $A_\infty$-module over $A$.

We define $\cc_\infty A$ to be the 
{\em category of $A_\infty$-modules}.
If $A$ is an ordinary associative algebra, we have a
faithful functor
\[
\cc A \ra \cc_\infty A
\]
but $\cc_\infty A$ has more objects and more morphisms
than $\cc A$. A typical example of an $A_\infty$-module
over  $A$ is obtained as follows: Suppose that
$A$ is the endomorphism algebra of a module $M'$
over an ordinary algebra $B$ and let $M$ be a
projective resolution of $M'$ over $B$. Then it is well-known
that each element of $A$ lifts to an endomorphism
of the complex $M$ and this endomorphism is unique up to homotopy.
The resulting `action up to homotopy' of $A$ on
$M$ may be shown to arise from a structure of
$A_\infty$-module, cf. \cite[6.1]{Keller99b}.

\bigskip
{\em 2nd step:} By definition, a morphism $f: L\ra M$ of
$A_\infty$-modules is {\em nullhomotopic} if there is a family
of graded maps
\[
h_n : L\ten A\tp{n-1} \ra M \ko n\geq 1\ko
\]
homogeneous of degree $-n$ such that
\begin{eqnarray*}
f_n
& = & \sum (-1)^{rs}\, m_{1+s} \circ (h_r \ten \id\tp{s})\\
&   & + \sum (-1)^{r+st}\, h_u \circ 
      (\id\tp{r}\ten m_s \ten \id\tp{t})\ko
\end{eqnarray*}
where the first sum runs over all decompositions
$n=r+s$, $r\geq 1$, $s\geq 0$ and the second sum over all
decompositions $n=r+s+t$, $r,t\geq 0$, $s\geq 1$ and
we put $u=r+1+t$.

The {\em homotopy category} $\ch_\infty A$ has
the same objects as $\cc_\infty A$ and morphisms
from $L$ to $M$ are morphisms of $A_\infty$-modules
modulo nullhomotopic morphisms.

\bigskip
{\em 3rd step:} To define the derived category, we should
formally invert all quasi-isomorphisms. This turns out
not to be necessary:

\begin{theorem} Each quasi-isomorphism of $A_\infty$-modules
is an homotopy equivalence.
\end{theorem}

\noindent
We therefore define the {\em derived category $\cd_\infty A$}
to be equal to the homotopy category $\ch_\infty A$. If $k$ is
only assumed to be a commutative ring, the assertion
of the theorem is no longer true in general and the derived category 
is different from the homotopy category.

\subsection{Solution of problem 1}

\label{solutionone}
\begin{theorem} Let $A$ be an $A_\infty$-algebra.
\begin{itemize}
\item[a)] If $A$ is an ordinary associative algebra with $1$,
the canonical functor 
\[
\cd A \ra \cd_\infty A
\]
is an equivalence onto the full subcategory of 
homologically unital $A_\infty$-modules, i.e. $A_{\infty}$-modules
$M$ such that $1\in A$ acts by the identity on $\Hs M$.
\item[b)] For each $A_\infty$-module $M$, the graded space
$\Hs M$ admits an $A_\infty$-module structure such
that $m_1=0$ and $M$ is isomorphic to $\Hs M$ in $\cd_\infty A$. 
It is unique up to (non unique, non strict) isomorphism of $A_\infty$-modules.
\end{itemize}
\end{theorem}

The $A_\infty$-module $\Hs M$ of b) is called the
{\em minimal model} of $M$. 
The theorem explains the answer to problem~1: Indeed, if
$M$ is a complex over an algebra $A$ as in (\ref{motivation}), 
then we may view it as an $A_\infty$-module, or, in other 
words, consider its image
in $\cd_\infty A$. By b), the image is isomorphic to
$\Hs M$ endowed with its canonical $A_\infty$-structure
$m_1=0, m_2, m_3, \ldots$. 
The full faithfulness of the functor $\cd A \ra \cd_\infty A$
then shows that $(\Hs M, m_1=0, m_2, \ldots)$ determines
the complex $M$ up to isomorphism in the derived category
$\cd A$.

In fact the theorem tells us more: 
If we also use the essential surjectivity
of the functor $\cd A \ra \cd_\infty A$ onto its image, we see that if
$V$ is a graded vector space, then we have a bijection between
complexes of $A$-modules $M$ with $\Hs M \iso V$ (as graded 
vector spaces) up to isomorphism in $\cd A$, and 
homologically unital $A_\infty$-structures $m_1=0, m_2, m_3, m_4, \ldots$ on 
$V$ up to isomorphism of such structures in $\cc_\infty\ca$.
Note that the latter set is naturally the set of orbits under the
action of a group, namely the group of families $f_i$, $i\geq 1$,
as in (\ref{MorphismAinfModule}) acting on the set of
homologically unital module structures by a conjugation action
deduced from (\ref{CompatMorphismAinfModule}). If
the total dimensions of $V$ and $A$ are finite, we even
obtain an algebraic group acting on an algebraic variety.

\section{Triangulated structure}
\label{triangstruct}

\subsection{Towards problem 2} \label{towardstwo}
Consider problem 2 of 
(\ref{motivation}) for the special case $n=1$: We then
have an ordinary associative algebra $B$ with $1$ and
a $B$-module $M$ and we wish to reconstruct the
category $\filt M\subset \Mod B$ of iterated extensions of 
$M$ by itself from the extension algebra 
\[
E=\Ext^*_B\,(M,M)
\]
with its canonical $A_\infty$-structure. The basic idea
is to realize $\filt M$ as a subcategory of the derived
category of $E$. In fact, this subcategory will simply be the 
subcategory of iterated extensions of the free $E$-module 
of rank one by itself. Now the derived category is 
not abelian in general so that a priori it makes no sense to speak
of extensions. However, like the derived category of an
ordinary associative algebra, the derived category 
of an $A_\infty$-algebra does admit a triangulated
structure and this suffices for our purposes.

\subsection{Triangulated structure on $\cd_\infty A$}
Recall that a {\em triangulated category} is an
additive category $\ct$ endowed with an additive
endofunctor $S:\ct \ra \ct$ called its {\em suspension}
(or {\em shift}) and with a class of sequences of the form
\[
X \ra Y \ra Z \ra SX\quad (*)
\]
called {\em triangles} (or {\em distinguished triangles}).
These have to satisfy a number of axioms 
(see the references at the beginning of \ref{reminder}).
If $Y$ occurs in a triangle $(*)$, it is called an
{\em extension} of $Z$ by $X$. 

For example, if $A$ is an ordinary associative algebra 
with $1$, the derived category $\cd A$ is triangulated:
The suspension of a complex $K$ is defined by
\[
(SK)^p=K^{p+1} \ko d_{SK}=-d_K.
\]
Each short exact sequence of complexes
\[
0 \ra K \ra L \ra M \ra 0
\]
gives rise to a canonical triangle
\[
K \ra L \ra M \ra SK
\]
and (up to isomorphism) all triangles are obtained
in this way. This generalizes to $A_\infty$-algebras
as follows

\begin{proposition} Let $A$ be an $A_\infty$-algebra.
Then the derived category $\cd_\infty A$ admits a
triangulated structure with suspension functor $S$
defined by
\[
(SM)^p = M^{p+1} \ko m_n^{SM}= (-1)^n\, m_n^M
\]
and such that each short exact sequence
\[
0 \ra K \arr{i} L \arr{p} M \ra 0
\]
of $A_\infty$-modules with {\em strict} morphisms
$i$ and $p$ gives rise to a canonical triangle.
Up to isomorphism in $\cd_\infty A$, all triangles
are obtained in this way.
\end{proposition}

\section{Standard functors}
\label{standardfunctors}

\subsection{More on problem 2} 
\label{moreonproblemtwo}
If $\ct$ is a triangulated
category and $M$ an object of $\ct$, we denote by
$\filt M=\filt_\ct M$ the full subcategory of $\ct$ whose objects
are iterated extensions of $M$ by itself. We denote
by $\tria M$ the closure of $M$ under suspension,
desuspension and forming extensions. Thus $\tria M$
is the smallest triangulated subcategory of $\ct$
containing $M$ and $\filt M$ is a full subcategory
of $\tria M$. 

If $\cs$ and $\ct$ are triangulated categories a
{\em triangle functor $\cs \ra \ct$} is an additive
functor $F: \cs \ra \ct$ together with an isomorphism
of functors $\phi: FS \ra SF$ such that for each triangle
$(X,u,Y,v,Z w)$ of $\cs$, the sequence
\[
FX \arr{Fu} FY \arr{Fv} FZ \longarr{(\phi X)(Fw)} SFX
\]
is a triangle of $\ct$. If $M$ is an object of $\ct$
and $F$ a triangle functor then clearly 
$F\tria M \subset \tria FM$ and similarly for $\filt M$.

In the situation of (\ref{towardstwo}), 
consider the subcategory $\filt E\subset \cd_\infty E$
obtained by taking the closure under extensions of the 
free $E$-module of rank one (also denoted by $E$).
We would like to show that the subcategory $\filt E$
is equivalent to $\filt_{\Mod B} M$. For this, we
first observe that the inclusion 
\[
\Mod B \ra \cd B
\]
induces an equivalence 
\[
\filt_{\Mod B} M \iso \filt_{\cd B} M.
\]
It will thus suffice to construct a triangle functor
$\cd B \ra \cd_\infty E$ inducing an equivalence
from $\filt M$ onto $\filt E$. The following diagram
summarizes the setup:
\[
\begin{array}{ccccc}
\Mod B & \longarr{ } & \cd B & \longarr{ } & \cd_\infty E \\
\ua          &             & \ua   &              & \ua    \\
\filt M     & \iso        & \filt_{\cd B} M & \iso & \filt E  
\end{array}
\]
The triangle functor $\cd B \ra \cd_\infty E$ is
constructed as a composition of standard functors.

\subsection{Restriction} \label{restriction}
Let $f: A \ra B$ be a morphism
of $A_\infty$-algebras and let $M$ be a module
over $B$. Define
\[
m^A_n : M\ten A\tp{n-1} \ra M\ko n\geq 1\ko
\]
by
\[
m_n^A=\sum (-1)^s\, m_{r+1}(\id\ten f_{i_1}\ten \cdots \ten f_{i_r})
\]
where the sum runs over all $1\leq r \leq n-1$ and all
decompositions $n-1=i_1+\cdots +i_r$ and the sign
is like in (\ref{morphisms}).
Then it is not hard to check that the $m_n^A$ define
an $A_\infty$-module $f^* M$ over $A$. The assignment
$M \ra f^*M$ becomes a functor
\[
\cc_\infty B \ra \cc_\infty A
\]
and induces a triangle functor
\[
\cd_\infty B \ra \cd_\infty A
\]
which will still be denoted by $f^*$. 

\begin{proposition} The functor $f^*:\cd_\infty A \ra \cd_\infty B$
is an equivalence if $f$ is a quasi-isomorphism.
\end{proposition}

\subsection{Hom and Tensor}
\label{HomAndTensor}
Let $B$ be an $A_\infty$-algebra and $X$ and $M$ 
two $B$-modules. For $n\in\Z$, let $\Hominf_B^n(X,M)$ be the space
of graded comodule maps $f: SX\ten TSB \ra SM \ten TSB$
of degree $n$ and define
\[
d(f)= b\circ f - (-1)^n f\circ b.
\]
Then $\Hominf_B^\bullet(X,M)$ is a complex which is bifunctorial
in $X$ and $M$. If $B$ is an ordinary algebra and $X,M$
are ordinary modules, then $\Hominf_B^\bullet(X,M)$
computes $\Ext^*_B\,(X,M)$.

Now suppose that $A$ is another 
$A_\infty$-algebra and that the $B$-module structure on $X$
comes from an $A$-$B$-bimodule structure, i.e. 
$X$ is endowed with graded maps
\[
m_{t,u}: A\tp{t-1}\ten X \ten B\tp{u-1} \ra X\ko t,u\geq 1\ko
\]
of degree $3-t-u$ such that we have the
associativity identity of (\ref{defainfalg}), where
the terms $\id\tp{r}\ten m_s \ten \id\tp{t}$ have
to be interpreted according to the type of their
arguments.
Then $\Hominf_B^\bullet(X,M)$ becomes a (right)
$A_\infty$-module over $A$ and we have the
Hom-functor
\[
\Hominf_B^\bullet(X, ?) :\cc_\infty B \ra \cc_\infty A
\]
\begin{proposition} The Hom-functor admits a left adjoint
\[
? \tensinf_A X : \cc_\infty A \ra \cc_\infty B.
\]
We have an induced pair of adjoint triangle functors
\[
\begin{array}{rcl} & \cd_\infty A  &\\
?\tensinf_A X & \da \ua & \Hominf_B^\bullet(X, ?) \\
 & \cd_\infty B &
\end{array}\quad .
\]
\end{proposition}

It is not hard to make $L\tensinf_A X$ explicit. Its
underlying graded space is $L\ten_k TSA \ten_k X$. If
$A$ is an ordinary associative algebra with $1$ and
$L,X$ are ordinary modules, then $L\tensinf_A X$ is 
a complex computing $\Tor_*^A\,(L,X)$.

\subsection{Application to problem 2} \label{application}
Let $B$ be an $A_\infty$-algebra and $M$ a $B$-module.
Then 
\[
A=\Hominf^\bullet_B\,(M,M)
\]
is a dg algebra and in particular an $A_\infty$-algebra.
The module $M$ becomes an $A$-$B$-bimodule.

\begin{proposition} The functor
\[
F=\Hominf^\bullet_B\,(M,?) : \cd_\infty B \ra \cd_\infty A
\]
induces a triangle equivalence
\[
\tria M \iso \tria A
\]
taking $M$ to $A$ and thus an equivalence $\filt M \iso \filt A$.
\end{proposition}

Now consider the situation of (\ref{towardstwo}). In this
case, $A$ is a dg algebra whose homology is the 
extension algebra of $M$. Let us simply write
\[
E=\Ext_B^*\,(M,M)
\]
for the extension algebra endowed with its canonical
structure of $A_\infty$-algebra. Then we know by Theorem
\ref{minimalmodels} that
we have a quasi-isomorphism $E \ra A$. This yields
the restriction functor
\[
\opname{res}: \cd_\infty A \ra \cd_\infty E
\]
which is an equivalence by Proposition \ref{restriction}. 
Consider the square
\[
\begin{array}{ccc} 
\cd_\infty B & \longarr{\sHominf^\bullet_B\,(M,?)} & \cd_\infty A \\
{\scriptsize\mathsf{can}}\ua \;\;&         &\da{\scriptsize \mathsf{res}} \\
\cd B  & \longarr{\quad\quad}                & \cd_\infty E
\end{array}
\]
The left vertical arrow induces an equivalence onto its 
image and maps $M$ to $M$. By the above proposition,
the top horizontal arrow induces an equivalence
$\filt M \iso \filt A$. 
The right vertical arrow is an equivalence mapping
the free module $A$ to the free module $E$ and hence
induces an equivalence $\filt A \iso \filt E$.
It follows that the bottom arrow induces an equivalence 
$\filt_{\cd B} M \iso \filt E$.
Finally, we know that $\filt_{\cd B} M$ is equivalent to
$\filt M\subset \Mod B$ so that we do obtain the required
equivalence between $\filt E \subset \cd_\infty E$ and
$\filt M \subset \Mod B$. 

\subsection{Solution of the original problem 2} \label{solution}
In order to treat problem 2 as it is stated in (\ref{motivation})
with $n$ modules $M_1, \ldots, M_n$, let $e_i$ be the idempotent
of $E$ corresponding to the direct factor $M_i$ of
$M=\bigoplus M_i$. There is a canonical isomorphism
\[
E^0 \iso \Hom_{\cd_\infty E}\,(E,E)
\]
and thus $e_i$ yields an idempotent endomorphism
of $E$ in the derived category. Now the derived category
is a triangulated category with infinite sums so that
by a well-known trick \cite{BoekstedtNeeman93}, the idempotent
$e_i$ splits in $\cd_{\infty} E$. Let $X_i$ denote its
image. We have 
\[
E \iso \bigoplus X_i
\]
in $\cd_\infty E$. The functor constructed in \ref{application}
sends $X_i$ to $M_i$ and thus induces an equivalence
$\filt X_i \iso \filt M_i$.

\section{Twisted objects}
\label{twistedobjects}

\subsection{Motivation}
Let $A$ be an associative algebra with $1$ and $\tria A$
the triangulated subcategory of $\cd A$ generated by the
free $A$-module of rank one (as in \ref{moreonproblemtwo}).
Up to isomorphism, the objects of $\tria A$ are the bounded
complexes of finitely generated free $A$-modules and
the morphisms between two such complexes are in bijection
with the homotopy classes of morphisms of complexes. 
In other words, $\tria A$ is {\em triangle equivalent
to the bounded homotopy category of finitely generated
free $A$-modules.} The advantage of this description
of $\tria A$ is that morphisms and objects are readily
computable. The aim of this section is to generalize
this description to any $A_\infty$-algebra
$A$ with a strict unit (cf. below). The generalization
will be of the form
\[
\tria A \iso \H{0} \tw A
\]
where $\tw A$ is the {\em $A_\infty$-category} of twisted
objects over $A$. An $A_\infty$-category is an `$A_\infty$-algebra
with several objects'  analogous to a `ring with several
objects' \cite{Mitchell72}. The above equivalence illustrates
the general fact that
working with several objects is natural and useful.
From the outset, we will therefore replace the $A_\infty$-algebra
$A$ by an $A_\infty$-category with strict identities $\ca$. 
We will then have more generally
\[
\tria \ca \iso \H{0}\tw \ca.
\]
The material in this section is adapted from
\cite{Kontsevich94}, \cite{Kontsevich98}. The special case
of dg categories goes back to Bondal-Kapranov's article
\cite{BondalKapranov90}.

\subsection{$A_\infty$-categories} 
\label{ainfcat}
Let $k$ be a field.
An {\em $A_\infty$-category $\ca$} is the datum of 
\begin{itemize}
\item a class of objects $\obj(\ca)$,
\item for all $A, B \in \ca$ (we write $A\in \ca$ instead 
of $A\in\obj(\ca)$), a $\Z$-graded vector space
$\Hom_\ca(A,B)$, often {\em denoted by $(A,B)$}, 
\item for all $n\geq 1$ and all $A_0, \ldots, A_n\in\ca$, a
graded map
\[
m_n : (A_{n-1}, A_n) \ten (A_{n-2},A_{n-1}) \ten
\cdots \ten (A_0, A_1) \ra (A_0, A_n)
\]
of degree $2-n$ 
\end{itemize}
such that for each $n\geq 1$ and all $A_0, \ldots, A_n\in\ca$,
we have the identity
\[
\sum (-1)^{r+st} \,m_u \,(\id\tp{r} \ten m_s \ten \id\tp{t}) = 0
\]
of maps
\[
(A_{n-1}, A_n) \ten (A_{n-2},A_{n-1}) \ten
\cdots \ten (A_0, A_1) \ra (A_0, A_n)\ko
\]
where the sum runs over all decompositions $n=r+s+t$ and
we put $u=r+1+t$. 

For example, the datum of an $A_\infty$-category $\ca$ with
one object $*$ is equivalent to the datum of the endomorphism
$A_\infty$-algebra $\Hom_\ca\,(*,*)$. Another class of
examples is provided by dg categories : these are exactly
the $A_\infty$-categories whose composition maps
$m_n$ vanish for all $n\geq 3$.

Note that for two reasons, an $A_\infty$-category is not
a category in general : firstly, the composition $m_2$ may
not be associative; secondly, there may not be
identity morphisms.

Let $\ca$ be an $A_\infty$-category and $A\in\ca$. A morphism
$e_A\in\Hom^0_\ca\,(A,A)$ is a {\em strict identity} if we
have
\[
m_2(f,e_A)=f\ko m_2(e_A,g)=g
\]
whenever these make sense and
\[
m_n(\ldots, e_A, \ldots) =0
\]
for all $n\neq 2$. In particular, $e_A$ is then a cycle
of the complex $(\Hom_\ca\,(A,A), m_1)$. Clearly, if $e_A$ exists,
it is unique. The $A_\infty$-category
{\em $\ca$ has strict identities} if there is a strict identity
$e_A$ for each object $A\in \ca$. 

\subsection{$A_\infty$-functors} 
Let $\ca,\cb$ be two $A_\infty$-categories. An {\em $A_\infty$-functor}
$F:\ca\ra\cb$ is the datum of
\begin{itemize}
\item a map $F:\obj(\ca) \ra \obj(\cb)$,
\item for all $n\geq 1$ and all $A_0, \ldots, A_n\in\ca$ a
graded map
\[
F_n : (A_{n-1}, A_n) \ten (A_{n-2},A_{n-1}) \ten
\cdots \ten (A_0, A_1) \ra \Hom_\cb\,(FA_0, FA_n)
\]
of degree $1-n$
\end{itemize}
such that conditions analogous to those of (\ref{morphisms}) hold.
The $A_\infty$-functor $F$ is {\em strict} if we have $F_n=0$ for
all $n\geq 2$. This implies that $F_1$ `commutes with all compositions'.

\subsection{Modules} 
We leave it as an exercise to the reader to generalize
(\ref{generalization}) to $A_\infty$-categories.
In particular, if $\ca$ is an $A_\infty$-category, he may
\begin{itemize}
\item define the category $\cc_\infty\ca$ of $\ca$-modules
\item define the derived category $\cd_\infty\ca$
\item define a structure of dg category
\[
\Hominf\,(L,M) \ko L,M\in \cc_\infty\ca
\]
on $\cc_\infty\ca$ as in \ref{HomAndTensor}
\item define the Yoneda functor
\[
Y : \ca \ra \cc_\infty \ca\ko A \mapsto \ca\,(?,A).
\]
Here we view $\cc_\infty\ca$ as a dg category (and hence
an $A_\infty$-category) and the Yoneda functor is to be
an $A_\infty$-functor.
\end{itemize}

\subsection{Factorization of the Yoneda functor} 
\label{factorization}
Let $\ca$ be an $A_\infty$-category and let $\tria \ca$
denote the triangulated subcategory of $\cd_\infty \ca$
generated by the image of the Yoneda functor.

\begin{theorem}
Suppose that $\ca$ has strict identities. Then there is
an $A_\infty$-category $\tw \ca$ and a factorization
%\[
%\begin{diagram}
%\node{\ca} \arrow{se,b}{Y}
%           \arrow{e,t}{Y_1} 
%\node{\tw\ca} \arrow{s,r}{Y_2} \\
%\node{ } \node{\cc_\infty \ca}
%\end{diagram} 
%\]
\[
\xymatrix{
\ca \ar[dr]_Y \ar[r]^{Y_1} &
\tw\ca \ar[d]^{Y_2} \\
&
\cc_\infty\ca
}
\]
of the Yoneda functor $Y$ such that $Y_1$
is a strict fully faithful $A_\infty$-functor and
$Y_2$ is an $A_\infty$-functor inducing an equivalence
\[
\H{0} \tw\ca \iso \tria\ca \subset \cd_\infty \ca.
\]
\end{theorem}
The theorem may be enhanced so as to {\em characterize
the category $\tw \ca$ of twisted objects by a
universal property}.

\subsection{Explicit description of $\tw \ca$} \label{explicit}
Let $\ca$ be an $A_\infty$-category with strict identities. 
We will explicitly construct the category of twisted objects $\tw \ca$
in two steps. First, we construct the $A_\infty$-category
$\Z\ca$, which is the `closure under shifts' of $\ca$:
the objects of $\Z\ca$ are the pairs $(A,n)$ of an
object $A\in\ca$ and an integer $n\in\Z$. The morphisms
from $(A,n)$ to $(A',n')$ are in bijection with
the elements of the graded space
\[
\Hom_\ca(A,A')[n'-n].
\]
We will identify $\ca$ with the full sub-$A_\infty$-category
of $\Z\ca$ on the objects $(A,0)$.
% COMPOSITION!
The compositions of $\Z\ca$ are obtained by shifting the
compositions of $\ca$. The reader is invited to write
down an explicit formula with correct signs.

In a second step, we construct $\tw \ca$ as
the `closure under extensions' of $\Z\ca$: Objects of $\tw \ca$ are
pairs $(B, \delta)$, where $B$ is a sequence $(A_1, \ldots, A_r)$
of objects of $\Z\ca$, and $\delta=(\delta_{ij})$ a matrix
of morphisms of degree $1$
\[
\delta_{ij} \in \Hom^1_{\Z\ca}\,(A_j,A_i)
\]
such that we have $\delta_{ij}=0$ for all $i\geq j$ and
\begin{equation}
\label{deftwistedobject}
\sum_{t=1}^\infty (-1)^{\frac{t(t-1)}{2}} 
m_t\,(\delta, \delta, \ldots, \delta)=0.
\end{equation}
Here we denote by $m_n$ the extension of $m_n^{\Z\ca}$ to
matrices. Note that since $\delta$ is strictly upper triangular,
the sum contains only finitely many non zero terms. 
The strange sign in the formula disappears in
the $b_n$-language (\ref{barconstruction}).

Alternatively, we can view an object $(B,\delta)$ as
the datum of an object $B$ of the additive hull
of $\Z\ca$ endowed with a finite split filtration 
and with an endomorphism $\delta$ of degree
$1$ which strictly decreases the filtration.

The space of morphisms from $(B,\delta)$ to $(B',\delta')$
is defined as
\[
\bigoplus_{i,j} \Hom_{\Z\ca}\,(A_i,A'_j).
\]
Let $n\geq 1$ and let $(B_0, \delta_0)$, \ldots, $(B_n,\delta_n)$
be objects of $\tw \ca$. We will define the composition map
(where the $\delta_i$ are suppressed from the notation)
\[
m_n: (B_{n-1}, B_n) \ten \cdots \ten (B_0, B_1) \ra (B_0, B_n).
\]
In the following formula, we write $\delta$ for all $\delta_i$
\[
m_n^{\tw\ca}= \sum_{t=0}^\infty \,\sum \, \pm m_{n+t}^{\Z\ca}\,\circ\,
(\id\tp{i_1} \ten \delta\tp{j_1} \ten \id\tp{i_2}\ten \delta\tp{j_2} \ten
\ldots \ten \id\tp{i_r} \ten \delta\tp{j_r}).
\]
Here the terms of the second sum are in bijection with
the non-commutative monomials
\[
X^{i_1} Y^{j_1} \ldots X^{i_r} Y^{j_r}
\]
of degree $n$ in $X$ and degree $t$ in $Y$ and the sign
is given by the identity
\[
(SX)^{i_1}(SY)^{j_1} \ldots (SX)^{i_r} (SY)^{j_r} =
\pm S^{n+t} X^{i_1} Y^{j_1} \ldots X^{i_r} Y^{j_r}
\]
in the algebra $k<X,Y,S>/(SX-XS, SY+YS)$. The  
pictorial version of the formula is
\[
m^{\tw\ca}_n=\sum \pm \centeps{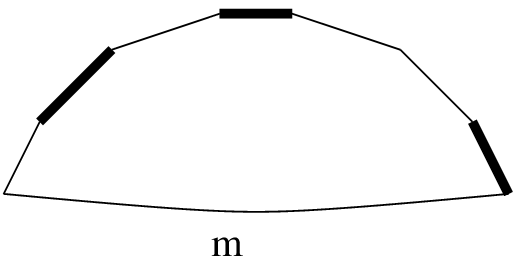} \;.
\]
Here the fat edges symbolize tensor powers of $\delta$'s
to be inserted in all possible places and in any number.
The concatenation of the thin edges represents a given
tensor $f_n\ten f_{n-1} \ten \cdots \ten f_1$.

It is worthwile to make the formula explicit 
for $m_1$. If $f$ is a morphism of even degree from
$(B,\delta)$ to $(B', \delta')$, we have
\begin{multline*}
m^{\tw \ca}_1(f) = m_1(f) - m_2(\delta', f) + m_2(f,\delta) \\
-m_3(\delta',\delta',f) + m_3(\delta', f,\delta) - m_3(f,\delta,\delta) \\
+m_4(\delta',\delta',\delta', f) - m_4(\delta',\delta',f,\delta)
+m_4(\delta', f, \delta,\delta) - m_4(f, \delta,\delta,\delta) + \ldots
\end{multline*}
If $f$ is of odd degree, we have 
\begin{multline*}
m^{\tw \ca}_1(f) = m_1(f)- m_2(\delta', f) - m_2(f,\delta) \\
-m_3(\delta',\delta',f) -m_3(\delta', f,\delta) - m_3(f,\delta,\delta) \\
+m_4(\delta',\delta',\delta', f) + m_4(\delta',\delta',f,\delta)
+m_4(\delta', f, \delta,\delta) + m_4(f, \delta,\delta,\delta) + \ldots
\end{multline*}
In the sequel, we identify $\ca$ with the subcategory 
of the objects $(A,0)$ of $\Z\ca$.
With this convention, the functor $Y_1$ of (\ref{factorization})
sends an object  $A\in\ca$ to $(A,\delta=0)$. The functor $Y_2$ sends
$B=(A_1, \ldots, A_n)$ to $YA_1\oplus \cdots \oplus YA_n$ 
(where $Y$ denotes the canonical extension of the Yoneda
functor to $\Z\ca$) endowed with the multiplications
\[
m_n = \sum_{t=0}^\infty (-1)^{\frac{t(t-1)}{2}} 
m_{n+t}\,\circ\,(\delta\tp{t} \ten \id\tp{n})\ko n\geq 1\ko
\]
where the $m_{n+t}$ on the right hand side denote the multiplications
of the direct sum of the $YA_i$.

\subsection{Description of filtered modules by twisted objects}
\label{DescriptionOfFilteredModules}
Let $B$ be an associative $k$-algebra
with $1$, let $M_1, \ldots, M_n$ be
$B$-modules and let $\filt(M_i)$ denote the full subcategory of
the category of right $B$-modules whose objects
admit finite filtrations with subquotients among
the $M_i$. Let $\ca$ be the $A_\infty$-category with
objects $1$, \ldots, $n$ and with morphism spaces
\[
\ca\,(i,j)= \Ext^*_{B}\,(M_i,M_j).
\]
The composition maps of $\ca$ are defined as follows:
The differential $m_1$ vanishes, the composition $m_2$ is the
Yoneda composition, and the $m_i$, $i>3$, 
are the canonical higher compositions constructed
as in (\ref{minimalmodels}, \ref{morphisms}).
We {\em assume} that these may be chosen so as
to make $\ca$ strictly unital. For example, this
holds if $\Ext^*_{B}\,(M_i,M_i)$ is one-dimensional
and $\Ext^*_{B}\,(M_i,M_j)=0$ for $i>j$ (cf. \ref{AugmentedAlgebras})
or if the $M_i$ are the simple modules over
a finite-dimensional basic algebra $B$. Indeed, in the latter
case, both $B$ and the dg algebra computing the extensions
are augmented over a product of fields (cf. \ref{AugmentedAlgebras}).

Let $\filt \ca$ denote the full subcategory of 
$\H{0}\tw \ca$ whose objects are the $(B,\delta)$
such that $B$ is a filtered direct sum of (non shifted)
copies of objects of $\ca$ (we use the
alternative description of $\tw\ca$ and identify $\ca$ with
a subcategory of $\Z\ca$).

\begin{theorem} There is an equivalence of
categories
\[
\filt\ca \iso \filt(M_i).
\]
\end{theorem}

Although the category $\filt\ca \subset \H{0}\tw \ca$
looks complicated, it is often quite
easy to describe in practice. 
We illustrate this in
two examples below.

The theorem is a consequence of (\ref{factorization}) and
(\ref{solution}). Indeed,
by its construction, the category $\ca$ is the 
$A_\infty$-$\Ext$-category of the $B$-modules
$M_1, \ldots, M_n$ in the sense of problem 2
of (\ref{motivation}). So by (\ref{application}) and
(\ref{solution}), there is a triangle functor
\[
F : \cd_\infty \ca \ra \cd B
\]
inducing an equivalence
\[
\tria(Y1, \ldots, Yn) \iso \tria(M_1, \ldots, M_n).
\]
By (\ref{factorization}), the Yoneda functor induces
an equivalence
\[
\H{0}\tw\ca \ra \tria(Y1,\ldots, Yn).
\]
Thus we get an equivalence
\[
\H{0}\tw\ca \ra \tria(M_1, \ldots, M_n)
\]
which takes the object $i$ to the corresponding
module $M_i$ and therefore induces an equivalence
\[
\filt \ca \iso \filt(M_i).
\]

\subsection{Example: Reconstruction from simple modules} 
\label{reconstrfromsimples} 
Let $B$ be a finite-dimensional basic
algebra. Theorem \ref{DescriptionOfFilteredModules} allows us
in particular to reconstruct the category $\mod B$ of
finite-dimensional $B$-modules from the
$A_\infty$-$\Ext$-category of the simple $B$-modules.

We illustrate this by an example. 
Let $B$ be the algebra given by the quiver with
relations 
\[
x \arr{\alpha} y \arr{\beta} z \arr{\gamma} t
\ko \gamma\,\beta\,\alpha=0
\]
and let $M_1=S_x, \ldots, M_4=S_t$ be the four simple $B$-modules.
The $A_\infty$-category $\ca$ of (\ref{DescriptionOfFilteredModules})
is then given by the quiver with relations
\[
\centeps{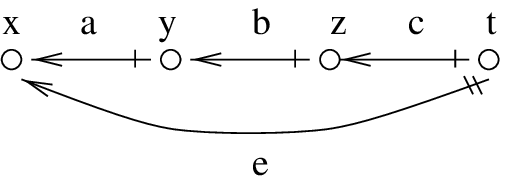} 
\]
where the arrows $a,b,c$ are of degree $1$, the arrow
$e$ is of degree $2$, and
all compositions other than $m_3(a,b,c)=e$ are defined
to be zero (this example has already appeared
in \ref{minimalmodels}). We would like to determine
the subcategory $\filt\ca$ of $\H{0}\tw\ca$
whose objects are the $(A,\delta)$ where $A$ is a
filtered sum of objects of $\ca$ (we use the
alternative description and identify $\ca$ with
a subcategory of $\Z\ca$). By an obvious abuse
of notation, we have
\[
A= V_x \ten x \oplus \ldots \oplus V_t\ten t
\]
where
$V_x, \ldots, V_t$ are finite-dimensional vector spaces
and the datum of
\[
\delta \in \Hom_{\Z\ca}^1\,(A,A)
\] 
corresponds to the datum of three linear maps $V_a$, $V_b$, $V_c$ :
\[
V_x \longlarr{V_a} V_y \longlarr{V_b} V_z \longlarr{V_c} V_t.
\]
The condition
\[
m_1(\delta) + m_2(\delta,\delta) - m_3(\delta,\delta,\delta) - \ldots =0
\]
translates into $V_a V_b V_c =0$. By definition, morphisms
from $(A,\delta)$ to $(A',\delta')$ are in bijection with
\[
\H{0} \Hom_{\tw \ca}\,((A,\delta), (A',\delta')).
\]
Now $\Hom_\ca^{-1}$ and hence 
\[
\Hom_{\tw \ca}^{-1}\,((A,\delta), (A',\delta'))
\]
vanishes
so that the homology classes in degree zero are exactly the zero
cycles. Now it is easy to see that the elements $f$ of
$\Hom^0_{\tw\ca}\,((A,\delta), (A',\delta'))$ correspond 
exactly to the quadruples
of linear maps
\[
f_x : V_x \ra V'_x \ko \ldots \ko f_t: V_t \ra V'_t
\]
and the condition
\[
m_1(f) - m_2(\delta', f) + m_2(f,\delta) - m_3( \ldots) + \ldots =0
\]
translates into
\[
f_x V_a = V'_a f_y \ko \ldots \ko f_z V_c = V'_c f_t.
\]
Hence, as expected, the category $\filt\ca$ is equivalent
to the category of contravariant representations of
the quiver with relations
\[
x \arr{\alpha} y \arr{\beta} z \arr{\gamma} t
\ko \gamma\,\beta\,\alpha=0\ko
\]
i.e.~to the category $\mod B$ of finite-dimensional
right modules over the algebra $B$. This confirms theorem
\ref{DescriptionOfFilteredModules} in this example.

\subsection{Example: Classifying filtered modules} \label{quasihered}
This example arose from an email exchange with T.~Br\"ustle.
Consider the following quiver $Q$
of type $D_4$ 
\[
\centeps{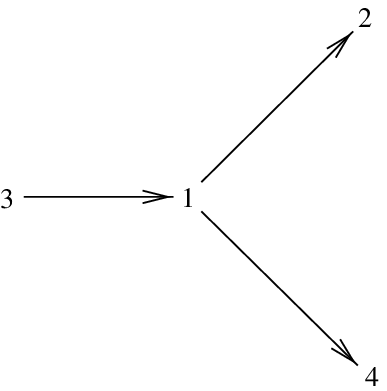}
\]
Denote by $S_i$ its simple (covariant) representation associated
with the point $i$, $1\leq i\leq 4$. Let $\Delta_1=S_1$,
$\Delta_2=S_2$, $\Delta_4=S_4$ and let $\Delta_3$ be the
following indecomposable representation 
\[
\centeps{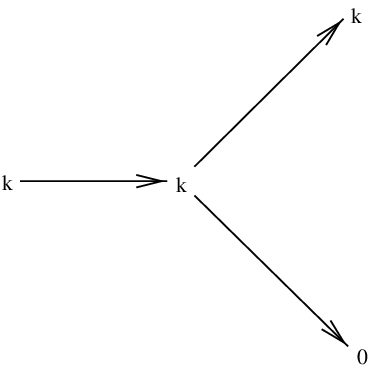} .
\]
We are interested in the category $\filt(\Delta_1, \ldots, \Delta_4)$
of representations of $Q$ admitting a filtration with 
subquotients among the $\Delta_i$. The $A_\infty$-$\Ext$-category
$\ca$ of the $\Delta_i$ is given by the quiver
\[
\centeps{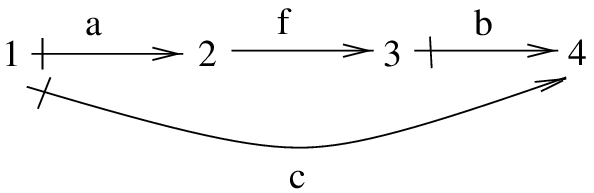}
\]
where $a,b,c$ are of degree $1$, $f$ is of
degree $0$ and all multiplications other than
$m_3(b,f,a)=c$ are defined to be zero. As in
example \ref{reconstrfromsimples}, the subcategory
$\filt\ca$ of $\H{0}\tw \ca$ can be
described as follows: Its objects are the finite-dimensional
(covariant) representations $V$ of the quiver
\[
2 \larr{a} 1 \arr{c} 4 \larr{b} 3
\]
of type $A_4$. Its morphisms $f:V\ra V'$ are the
quintuples of linear maps $f_i: V_i \ra V'_i$
and $\phi: V_2 \ra V'_3$ 
%\[
%\begin{diagram}
%\node{V_2} \arrow{s,l}{f_2} \arrow{esese,t}{\phi}
%\node{V_1} \arrow{w,t}{V_a} \arrow{s,r}{f_1} \arrow{e,t}{V_c} 
%\node{V_4} \arrow{s,r}{f_4} 
%\node{V_3} \arrow{w,t}{V_b} \arrow{s,r}{f_3} \\ \\`
%\node{V'_2} 
%\node{V'_1} \arrow{w,b}{V'_a} \arrow{e,b}{V'_c} 
%\node{V'_4} 
%\node{V_3} \arrow{w,b}{V'_b}
%\end{diagram}
%\]
\[
\xymatrix{V_2  \ar[d]_{f_2} \ar[drrr]^{\phi} &
V_1 \ar[l]_{V_a} \ar[d]_{f_1} \ar[r]^{V_c} &
V_4 \ar[d]^{f_4} &
V_3 \ar[d]^{f_3} \ar[l]^{V_b}\\
V'_2  &
V'_1 \ar[l]^{V'_a} \ar[r]_{V'_c} &
V'_4 &
V'_3 \ar[l]^{V'_b}
}
\]
satisfying
\begin{align*}
0= V'_a \,f_1- f_2 V_a \ko 0 = V'_b f_3 - f_4 V_b \\
0 = V'_c f_1 - f_4 V_c + V'_b \,\phi\, V_a \quad (*) 
\end{align*}
Thanks to theorem \ref{DescriptionOfFilteredModules}, we have an
equivalence
\[
\filt\ca\iso \filt(\Delta_1, \ldots, \Delta_4).
\]
The term $V'_b\,\phi\,V_a$ in $(*)$ is essential: If it
is omitted, the problem is that of classifying
the representations of the quiver of type $A_4$ above.
Then we find 10 isoclasses of indecomposables
instead of the 9 isoclasses present on
the right hand side. 

\section{On the proofs}

\subsection{The minimality theorem \ref{minimality}}
The article \cite{JohanssonLambe96} gives a detailed account of
various approaches to this theorem (for the case of connected dg
algebras $A$).  For the case where $A$ is a $\Z$-graded dg algebra, it
is proved in \cite{Merkulov98}. Here one also finds new explicit
formulas for the higher multiplication maps on $\Hs A$.

In the general case, one may proceed as follows \cite{Lefevre2000}:
Let $A$ and $B$ be two $A_\infty$-algebras and $(f_1, \ldots, f_s)$ an
$A_s$-morphism (i.e. the defining equation from \ref{morphisms} holds
for $n\leq s$). One can then define a morphism of complexes $c=c(f_1,
\ldots, f_s): A\tp{s+1} \ra B$ such that $(f_1,\ldots,f_s)$ extends to
an $A_{s+1}$-morphism $(f_1, \ldots, f_{s+1})$ iff $c$ is null
homotopic and in this case, the homotopy is given by $f_{s+1}$.  Here,
$A$ and $B$ are endowed with the differentials given by $m_1$. It
follows that if $(B,m_1)$ is contractible, each $A_s$-morphism with
target $B$ extends to an $A_\infty$-morphism.

Now suppose that $A$ is an $A_\infty$-algebra and $A=A'\oplus N$ a
decomposition of the complex $(A,m_1)$ with contractible $N$.  We wish
to make $A'$ into an $A_\infty$-algebra quasi-isomorphic to $A$. For
this, we make $N$ into an $A_\infty$-algebra by taking all $m_n$,
$n\geq 2$, to vanish.  Since $N$ is contractible, the morphism of
complexes $A \ra N$ extends to an $A_\infty$-morphism. We thus obtain
a dg coalgebra morphism $p : \ol{T}SA \ra \ol{T}SN$ (cf.
\ref{barconstruction}).  It is easy to check that its kernel in the
category of graded coalgebras is isomorphic to $\ol{T}SA'$. 
Thanks to \ref{barconstruction}, this
yields an $A_\infty$-algebra structure on $A'$ and an
$A_\infty$-morphism $A' \ra A$ extending the inclusion.

It is easy to see that an $A_\infty$-morphism $f$ between
minimal $A_\infty$-algebras is invertible iff $f_1$ is
invertible. This yields the uniqueness assertion of the theorem.

\subsection{The theorem on homotopies \ref{homotopyequivalences}}
The theorem results from the facts, proved in \cite{Lefevre2000},
that 
\begin{itemize}
\item[1)] the bar construction functor (\ref{barconstruction}) is
an equivalence from the category of $A_\infty$-algebras onto
the full subcategory of 
fibrant-cofibrant objects of a closed model category of coalgebras and 
\item[2)] the intrinsic homotopy relation of this 
closed model category induces the explicit relation given
in \ref{homotopyequivalences}.
\end{itemize}

Alternatively, one can prove b) by the following argument taken
from \cite[4.5]{Kontsevich97}:
One shows that the quasi-isomorphism $i_A : \Hs A \ra A$ constructed
above is actually an homotopy equivalence. The given morphism
$s: A \ra B$ thus fits into a diagram
\[
\xymatrix{
A \ar[r]^s  & B \\
\Hs A \ar[u]^{i_A} \ar[r]_{s'} & \Hs B \ar[u]_{i_B}
}
\]
of $A_\infty$-algebras which is commutative up to homotopy.
The first component $s'_1$ identifies with $\Hs s$. So
$\Hs s$ is invertible iff $s'$ is an isomorhism iff
$s$ is a homotopy equivalence.

\subsection{$A_\infty$-modules  (section \ref{ainfmodules})} 
Theorem \ref{generalization} on $A_\infty$-modules
has an analogous (but simpler) proof as the corresponding 
theorem \ref{homotopyequivalences}
on $A_\infty$-algebras, cf. \cite{Lefevre2000}.

Part a) of theorem \ref{solution} is proved in \cite[section 7]{Keller99b}.
The existence of the minimal model in b) is shown
in \cite{GugenheimLambeStasheff91}. Uniqueness is easy
(as in the case of $A_\infty$-algebras).

\subsection{Triangulated structure (section \ref{triangstruct})} 
The following argument is
taken from \cite{Lefevre2000}: The category of $A_\infty$-modules
over an $A_\infty$-algebra $A$ is equivalent to a full subcategory
$\cc'$ of the category $\cc$ of differential counital comodules over
the counital dg coalgebra $TSA$ (cf. \ref{generalization}).  We endow
the category $\cc$ with the class of exact sequences whose underlying
sequences of graded comodules split.  As in \cite[2.2]{Keller94}, one
shows that $\cc$ is a Frobenius category. The associated stable
category is thus triangulated.  This category equals the homotopy
category $\ch$ of the category of dg counital comodules. It is easy to
check that the essential image $\ch'$ of $\cc'$ in $\ch$ is a
triangulated subcategory. By theorem \ref{generalization}, $\ch'$ is
equivalent to the derived category $\cd_\infty A$.  Exact sequences of
strict morphisms of $A_\infty$-modules give rise to exact sequences of
coinduced graded $TSA$-comodules.  Such sequences split in the
graded category (coinduced modules are injective). Thus each
exact sequence of strict morphisms of $A_\infty$-modules gives rise to
a canonical triangle. When we have an arbitrary triangle $(K,L,M)$ of
$\ch'$, we can form a mapping cylinder over the morphism $M\ra SK$ to
obtain an isomorphic triangle $K \ra L' \ra M \ra SK$.  By the
construction of the mapping cylinder, the sequence of graded comodules
$K\ra L' \ra M$ identifies with the canonical sequence $K \ra K\oplus
M \ra M$. So the sequence $K \ra L' \ra M$ does come from a sequence
of strict $A_\infty$-morphisms.

\subsection{Standard functors (section \ref{standardfunctors})} 
The assertions in this section are based on the following description
\cite{Lefevre2000} of the derived category $\cd_\infty A$ :
The functor $M\mapsto M\ten TSA$ is an
equivalence of the category of $A_\infty$-modules onto the
full subcategory of fibrant-cofibrant objects of a closed
model category of counital dg comodules over $TSA$. Thus
it induces an equivalence of $\cd_\infty A$ onto the 
corresponding homotopy category. 

\subsection{Twisted objects (section \ref{twistedobjects})}
We follow \cite{Lefevre2000}.
Below, we will prove that $\tw\ca$ is an $A_\infty$-category.
The proof that $Y_2$ is an $A_\infty$-functor follows the
same lines. It remains to be proved that $Y_2$ induces an equivalence
of $H^0 \tw\ca$ onto its image in the homotopy category
(=derived category) of $A_\infty$-modules. For this, let $B_1$ and $B_2$ be two
twisted objects and consider the map
\[
H^* \tw\ca(B_1, B_2) \ra H^* \Hominf_\ca(Y_2 B_1, Y_2 B_2)
\]
induced by $Y_2$. We have to show it is an isomorphism.
Using the canonical filtrations of $B_1$
and $B_2$, one easily reduces this to the special case
where the $B_i$ are objects $A_i$ of $\ca$ (and thus the
corresponding twisting elements $\delta$ vanish). Then one
has to show that the map
\[
H^* \ca(A_1, A_2) \ra H^* \Hominf_\ca(YA_1, YA_2)
\]
is bijective. Now using the strict identity of $A_1$, it is
not hard to construct an explicit inverse of this map.

Let us show that $\tw\ca$ is an $A_\infty$-category.
Let $\cb$ be the datum of
\begin{itemize}
\item the same class of objects as $\tw\ca$,
\item the same graded morphism spaces as $\tw\ca$,
\item the natural composition maps $m'_n$ obtained by extending
the maps $m_n^{\Z\ca}$ to matrices.
\end{itemize}
It is clear that $\cb$ is an 
$A_\infty$-category. For each sequence $B_0, \ldots, B_n$ of
objects of $\cb$, let $b'_n$ be the map
\[
S(B_{n-1},B_n) \ten S(B_{n-2},B_{n-1}) \ten \ldots \ten S(B_0,B_1)
\ra S(B_0,B_n)
\]
deduced from $m'_n$ as in the bar construction
(\ref{barconstruction}). Here, $S(B_{i-1},B_i)$ denotes the suspension
of the morphism space $(B_{i-1},B_i)$. For each $0\leq i\leq n$, let
$\alpha : k \ra S(B_i,B_i)$ be the composition of the degree $1$ map
\[
k \ra (B_i, B_i)\ko 1\mapsto \delta
\]
with the canonical morphism $s:(B_i,B_i) \ra S(B_i,B_i)$
of degree $-1$. Note that $\alpha$ is of degree $0$. 
A sign check shows that the identity \ref{deftwistedobject} translates
into
\[
0 = \sum_{i\geq 1} b'_i\,\circ \,\alpha\tp{i}\, .
\]
Now for each $t\geq 0$, let $\phi_{n,t}$ be the sum of the maps
\[
\id\tp{i_1} \ten \alpha\tp{j_1} \ten \ldots \ten
\id\tp{i_r}\ten\alpha\tp{j_r}
\]
defined on 
\[
S(B_{n-1},B_n) \ten \ldots \ten S(B_0,B_1)
\]
where $i_1+\cdots+i_r=n$ and $j_1+\cdots+j_r=t$ as in 
\ref{explicit}. 
Note that for each decomposition
$n=n_1+n_2$, we have
\begin{equation} \label{decompositionformula}
\phi_{n,t} = \sum_{t=t_1+t_2} \phi_{n_1,t_1} \ten \phi_{n_2,t_2}.
\end{equation}
Define
\[
b_n = \sum_{t\geq 0} b'_{n+t}\circ \phi_{n,t}.
\]
A sign check shows that $b_n$ corresponds
to the twisted composition $m_n^{\tw\ca}$ 
defined in \ref{explicit}. We have to show that
\[
\sum b_l\, (\id\tp{i} \ten b_j \ten \id\tp{k'}) =0  \ko
\]
where the sum runs over the decompositions $n=i+j+k'$ and
we put $l=i+1+k'$. Now we have
\begin{align*}
\sum b_l\, (\id\tp{i} \ten b_j \ten \id\tp{k'})
& =
\sum
b'_{l+s} \, \phi_{l,s}\, (\id\tp{i} \ten (b'_{j+t} \,\phi_{j,t})\ten
\id\tp{k'}) \\
& =
\sum
b'_{l+s}\, (\phi_{i,s_1} \ten (b'_{j+t} \,\phi_{j,t}) \ten
\phi_{k',s_2})
\end{align*}
the indices run trough the integers $\geq 0$ such that $s=s_1+s_2$.
Here, we have used the decomposition formula
\ref{decompositionformula} and the fact that $\alpha$ is of degree
$0$.  Let $\Sigma_1$ denote the last sum.  Graphically, the terms of
$\Sigma_1$ correspond to diagrams
\[
\centeps{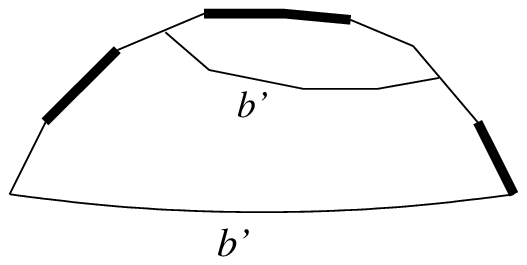} \;,
\]
where the fat edges symbolize $\alpha$'s to be inserted in any
number and in all positions such that the arc above the upper $b'$ 
is not entirely fat (it must always contain a thin piece 
symbolizing the $j\geq 1$ arguments of $b_j$).
To show that $\Sigma_1$ vanishes, we compare it to
\begin{equation}
\sum_{r\geq 0} \sum b'_L\,(\id\tp{I} \ten b'_J \ten \ko
\id\tp{K})\,\phi_{n,r}
\end{equation}
where the inner sum runs over all decompositions of $N=n+r$
into $N=I+J+K$ and we put $L=I+1+K$. 
Since $\cb$ is an $A_\infty$-category,
the inner sum vanishes for each $r$.
On the other hand, the whole expression equals
\[
\sum b'_L\,(\phi_{I', r_1} \ten (b'_{J'+r_2}\, \phi_{J', r_2}) 
\ten \phi_{K',r_3}) \ko
\]
where the sum runs over the non negative integers
$r, r_1, r_2, r_3, I', J', K'$ satisfying
$r=r_1+r_2+r_3$, $I'+J'+K'=n$, and we put $L=(I'+r_1)+1+(K'+r_3)$.
Let us denote the last sum by $\Sigma_2$. 
Each term of $\Sigma_1$ appears exactly once in $\Sigma_2$
but not all terms of $\Sigma_2$ appear in $\Sigma_1$. Those
that do not appear are precisely the terms where $J'=0$, i.e. the
terms
\[
b'_L\, (\phi_{I',r_1} \ten (b'_{r_2}\, \phi_{0, r_2}) \ten
\phi_{K',r_3}).
\]
Graphically, they correspond to diagrams
\[
\centeps{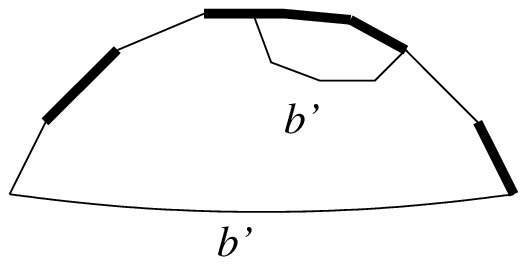} \;.
\]
Now we have $b'_{r_2}\,\phi_{0,r_2} = b'_{r_2}\,\circ\,\alpha\tp{r_2}$
so that the sum over all $r_2\geq 0$ of these terms vanishes.
Thus we have $\Sigma_1=\Sigma_2=0$.

\end{document}